\documentstyle{amsart}
\theoremstyle{plain}
\newtheorem{theorem}{Theorem}
\newtheorem{lemma}[theorem]{Lemma}
\newtheorem{proposition}[theorem]{Proposition}

\newtheorem{corollary}[theorem]{Corollary} 
\numberwithin{equation}{section}
\numberwithin{theorem}{section}

\theoremstyle{definition}

\theoremstyle{remark}

\begin{document}

\title[Zero-pole interpolation]{Zero-pole interpolation for matrix 
meromorphic functions on a 
compact Riemann surface and a matrix Fay trisecant identity} 
\author{Joseph A. Ball}
\address{Department of Mathematics \\ Virginia Tech \\
Blacksburg, Virginia 24061}
\email{ball@@math.vt.edu}
\thanks{The first author was supported by the NSF under grant DMS-9500912 
and as a member of the ``Holomorphic Spaces'' program sponsored 
by the Mathematical Sciences Research Institute in Fall 1995.}
\author{Victor Vinnikov}
\address{Department of Theoretical Mathematics\\ 
Weizmann Institute of Science\\
Rehovot 76100, Israel}
\email{vinnikov@@wisdom.weizmann.ac.il}
\thanks{The second author was also supported as a member of the 
``Holomorphic Spaces'' program sponsored by the 
Mathematical Sciences Research Institute in Fall 1995,
 and is an incumbent of Graham 
and Rhona Beck career development chair, Weizmann Institute.}
\keywords{flat vector bundles, connection, Cauchy kernel, factor of 
automorphy}
\subjclass{Primary: ; Secondary: }

\maketitle

ABSTRACT:  This paper presents a new approach to constructing a 
meromorphic bundle map between flat vector bundles over a compact Riemann 
surface having a prescribed Weil divisor (i.e. having prescribed zeros 
and poles with directional as well as multiplicity information included 
in the vector case).  This new formalism unifies the earlier approach of 
Ball-Clancey (in the setting of trivial bundles over an abstract Riemann 
surface) with an earlier approach of the authors (where the Riemann 
surface was assumed to be the normalizing Riemann surface for an 
algebraic curve embedded in ${\bold C}^2$ with determinantal 
representation, and the vector bundles were assumed to be presented as 
the kernels of linear matrix pencils).  The main tool is a version of the 
Cauchy kernel appropriate for flat vector bundles over the Riemann 
surface. Our formula for the interpolating bundle map (in the special 
case of a single zero and a single pole) can be viewed as a 
generalization of the Fay trisecant identity from the usual 
line bundle case to the vector bundle case in terms of Cauchy kernels.
In particular we obtain a new proof of the Fay trisecant identity.

\section{Introduction} \label{S:intro}

The following zero-pole interpolation problem is one of the main objects 
of study in the recent monograph \cite{bgr}.  We state here the simplest 
case where all zeros and poles are assumed to be simple and disjoint.  
{\sl Given a finite collection $\lambda^1, \dots, 
\lambda^{n_0}, \mu^1, \dots, \mu^{n_\infty}$ of distinct points in the 
complex plane ${\bold C}$, nonzero column vectors $u_1, \dots, 
u_{n_\infty} \in {\bold C}^{r \times 1}$ and nonzero row vectors $x_1, 
\dots, x_{n_0} \in {\bold C}^{1 \times r}$, find (if possible) an $r 
\times r$ matrix function $T(z)$ having value equal to the identity matrix
$I$ at infinity such that 
(1) $T(z)$ is analytic on $({\bold C} \cup \{\infty\}) 
\backslash \{\mu^1, \dots, 
\mu^{n_\infty}\}$ and $T(z)^{-1}$ has analytic continuation to $({\bold 
C} \cup \{\infty \}) \backslash \{\lambda^1, \dots, \lambda^{n_0} \}$, 
(2) for $i=1, 
\dots, n_0$, $T(\lambda^i)$ has rank $r-1$ and $x_i T(\lambda^i) = 0$, 
and (3) for $j=1, \dots, n_\infty$, $T(\mu_j)^{-1}$ (i.e., the analytic 
continuation of $T(z)^{-1}$ to $z=\mu^j$) has rank $r-1$ and $T(\mu^j)^{-1} 
u_j = 0$.}  A complete solution, along with numerous applications to 
problems in factorization, matrix interpolation and $H_\infty$-control, 
is given in \cite{bgr}.  The solution (for this simple case) is 
as follows.  {\sl A solution exists if and only if the $n_0 \times 
n_\infty$ matrix 
\begin{equation}
\Gamma = [\Gamma_{ij}]\text{ with } \Gamma_{ij} = \dfrac{x_i u_j}{\mu^j - 
\lambda^i}
\label{Gamma}
\end{equation}
 is square and invertible.  In this 
case the unique solution is given by
\begin{equation}
T(z) = I + \sum_{j=1}^{n_\infty} u_j (z-\mu^j)^{-1} c_j
\label{genus0solution}
\end{equation}
where $c= \begin{bmatrix} c_1 & \dots & c_{n_\infty} \end{bmatrix} ^T$ is 
given by}
\begin{equation}
c= \Gamma^{-1} \begin{bmatrix} x_1^T & \dots & x_{n_0}^T \end{bmatrix}^T.
\label{scalarcoefficients}
\end{equation}

The solution in \cite{bgr} uses system theory ideas, especially the state 
space similarity theorem specifying the level of uniqueness for two 
realizations of the same rational matrix function as the transfer function 
of a linear system.  Later work (see \cite{bgrak}) handles the nonregular 
case (where $\det T(z)$ vanishes identically and the nature of the zero 
structure must be enlarged) by elementary linear algebra, without 
recourse to the state space similarity theorem. 

There have now appeared two seemingly distinct generalizations of this 
result to higher genus.  In \cite{bc}, the problem is posed to construct 
a (global, single-valued) meromorphic matrix function on the compact 
Riemann surface $X$ satisfying conditions as in (1), (2) and (3) above.  
A matrix analogous to $\Gamma$ appears, but the solution criterion is not 
as simple; nevertheless, an explicit formula analogous to 
\eqref{genus0solution} and \eqref{scalarcoefficients} was found for the 
solution when it exists.  The approach in \cite{bc} can be seen as an analogue 
of that in \cite{bgrak} (i.e., system theory ideas are avoided
and a simple ansatz is used to reduce the problem to an analysis of a 
linear system of equations).  The 
paper \cite{hip}, on the other hand, while formulating a more general 
problem (involving bundle maps between certain types of flat vector 
bundles rather than global meromorphic matrix functions) in an abstract 
setting, works primarily in a more concrete setting, where the Riemann 
surface $X$ is taken to be the normalizing Riemann surface for an 
algebraic curve $C$ having a determinantal representation
\[
C = \{ (\lambda_1, \lambda_2) \in {\bold C}^2 \colon \det (\lambda_1 
\sigma_2 - \lambda_2 \sigma_1 + \gamma) = 0 \}
\]
(where $\sigma_1, \sigma_2, \gamma$ are $M \times M$ matrices and 
$\lambda = (\lambda_1, \lambda_2)$ are affine coordinates), 
and the input and output 
bundles $E$ and $\widetilde{E}$ are assumed to have kernel 
representations, e.g.,
\[
E(\lambda) = \{ v \in {\bold C}^M \colon (\lambda_1 \sigma_2 - \lambda_2 
\sigma_1 + \gamma) v = 0\}.
\]  
In this 
setting, a non-metric version of the several-variable system theory 
connected with the model theory for commuting operators due to Livsic 
(i.e., a version with all Hilbert space inner products dropped) applies, 
any meromorphic bundle map satisfying appropriate conditions at the 
points at infinity can be realized as the transfer function, or the 
joint characteristic function, of a Livsic-Kravitsky 2D system, and the 
zero-pole bundle-map interpolation problem can be solved using the 
state-space similarity theorem for this setting in a manner 
completely parallel to that of 
\cite{bgr} for the genus 0 case.  (For a recent systematic treatment of 
the Livsic theory, we refer to \cite{lkmv} and \cite{vinsurvey}).  
In this solution there is a 
matrix $\Gamma$ analogous to the $\Gamma$ in \eqref{Gamma} along with an
explicit formula for the solution (when it exists) as in 
\eqref{genus0solution}. In this setting, one specifies an output bundle 
$\widetilde{E}$ having a kernel bundle representation as well as the zero-pole 
interpolation data.  The invertibility of $\Gamma$ is then equivalent to 
the existence of an input bundle $E$ also having a kernel representation 
together with a bundle map $T \colon E \to \widetilde{E}$ meeting the 
zero-pole interpolation conditions.

The purpose of this paper is to synthesize these two approaches.  We 
obtain a generalization of the approach of \cite{bc} which handles the 
vector bundle problem, and clarify the solution criterion as well as the 
role of the invertibility of $\Gamma$ in this abstract setting.  To 
obtain an analogue of the basic ansatz in \cite{bc} used for the form of 
the solution for the general bundle case,  we need a version of the 
Cauchy kernel $(z,w) \to \frac{1}{z-w}$ for sections of a flat vector 
bundle $\chi$ satisfying $h^0(\chi \otimes \Delta) = 0$ where $\Delta$
is a line bundle of differentials of order 1/2 (a theta characteristic or 
a spin structure).  
This object was introduced in \cite{hip} (see also \cite{AV2}, \cite{vin1}
for the line bundle case) 
but a proof of its existence for the vector bundle case relied on the 
theory of determinantal representations of algebraic curves and of kernel 
representations of bundles over such curves.  Here we give a simple, 
direct existence proof using only some cohomology theory of vector 
bundles and the Riemann-Roch theorem.  A similar proof of the same result 
for the line bundle case is given in \cite{raina1} and \cite{raina2}; the 
general case is also handled in \cite{newfay} but by using completely 
different techniques (involving the theory of the Green's function for 
the heat equation over $X$).  Various other forms of the Cauchy kernel 
for a Riemann surface have appeared earlier in the literature, in 
particular in connection with the Riemann-Hilbert problem (see 
\cite{Rodin}, \cite{Zverovich}).  However, these are developed within the 
framework of meromorphic differentials whereas our Cauchy kernel is 
defined as a multiplicative meromorphic differential of order $1/2$.  The 
use of half-order differentials has the advantage that no extraneous 
poles are introduced in the Cauchy kernel.  We mention that the paper 
\cite{AV} applies our Cauchy kernel to the study of indefinite Hardy 
spaces on a finite bordered Riemann surface.

Taking the constant term in the Laurent expansion of the Cauchy kernel 
around the diagonal allows us to define a certain flat connection on the 
flat bundle $\chi$.  In the concrete setting of the determinantal 
representations, this connection was already introduced in \cite{hip}.  
This flat connection is determined canonically up to a choice of a bundle 
$\Delta$ of 
half-order differentials.

We also make explicit the mappings between the concrete and abstract 
settings; in this way we are able to see explicitly the equivalence 
between the solution in \cite{hip} and the solution in \cite{bc}.
The main ingredient is the explicit formula for the determinantal 
representation of an algebraic curve 
with a given kernel bundle in terms of the Cauchy kernel 
of the bundle.

We also specialize the results to the line bundle case.  For this case, 
both the solution to the zero-pole interpolation problem and the Cauchy 
kernel can be expressed explicitly in terms of theta functions (see 
\cite{oldfay}, \cite{mumford} and \cite{farkaskra} for background 
material on theta functions).  When this is done, the equality between 
these two forms of the solution of the interpolation problem leads to a 
new proof of the trisecant identity due to Fay (see \cite{oldfay} 
Corollary 2.19 or \cite{mumford} Volume II page 3.214). 
In the general vector bundle case, the formula for the solution of the 
interpolation problem in terms of the Cauchy kernels (in the case of a 
single zero and a single pole---see \eqref{3.12a}) can be viewed as a 
matrix version of the Fay trisecant identity.
 We mention that 
in the genus 1 case one can obtain an explicit formula for the Cauchy 
kernel for the general case of flat vector bundles (see \cite{bcv}).

We close the introduction by mentioning three other possible further 
applications of our Cauchy kernel.

First of all, as will be shown in Section \ref{S:detrep}, the vector 
bundle $\chi$ is completely determined by the values $K(\chi; x^i,x^j)$ 
of the Cauchy kernel at a certain finite collection of points $x^1, 
\dots, x^m$ (forming a line section in a birational planar embedding of 
$X$).  This suggests that these values can be used as affine coordinates 
for the bundle $\chi$ in the corresponding 
moduli space of semistable bundles on the complement of the generalized 
theta divisor.  A very 
similar construction for line bundles on hyperelliptic curves is due to 
Jacobi (see \cite{Jacobi}) and has been given a modern treatment by 
Mumford (see Volume II of \cite{mumford}).

Secondly, in the line bundle case consideration of the Fay trisecant 
identity when some of the points come together leads to very interesting 
identities, showing in particular that the theta function satisfies the 
KP equations (see \cite{oldfay}, Volume II of \cite{mumford} and 
\cite{Shiota}).  An interesting line of research is to consider similar 
limiting versions of our matrix Fay trisecant identity \eqref{3.12a}.
A related problem is to find the relation between the Cauchy kernel
and the matrix Baker-Akhiezer function of Krichever and Novikov \cite{KrNo}
whose definition involves the so called Tjurin parameters \cite{Tju1,Tju2}
of the vector bundle.

Thirdly, the absence of explicit formulas for the Cauchy kernel makes it 
interesting to try to find formulas for the Cauchy kernel of one vector 
bundle in terms of the Cauchy kernel for another.  In particular it would 
be interesting to find how the Cauchy kernel behaves under pullback and 
direct image.

The paper is organized as follows.  Section \ref{S:intro} is this 
introduction.  Section \ref{S:cauchyker} develops the Cauchy kernel for a 
flat line bundle.  Section \ref{S:absint} then formulates and solves the 
zero-pole interpolation problem in the abstract setting.  Section 
\ref{S:linebundle} obtains explicit formulas for all the results in the 
line bundle case and obtains the new proof of the Fay trisecant 
identity.   Section \ref{S:detrep} explains how to use Cauchy kernels to 
obtain a canonical map from the abstract to the concrete setting.  
Finally, Section \ref{S:conint} explains the connections with the 
concrete interpolation problem solved in \cite{hip}.

\section{The Cauchy kernel for a flat vector bundle} 
\label{S:cauchyker}

We assume that we are given a compact Riemann surface $X$. 
Let $\Delta$ be a line bundle of differentials of order 
$\frac{1}{2}$ on $X$, i.e., a line bundle satisfying $\Delta \otimes 
\Delta \cong K$, where $K$ is the canonical line bundle (i.e., the 
line bundle with local holomorphic sections equal to local holomorphic 
differentials on $X$).  Note that since $\deg(K) = 2g-2$, $\deg(\Delta) = 
g-1$ where $g$ is the genus of $X$.
In addition we assume that we are given a holomorphic
complex vector bundle $\chi$ of degree $0$ (and of rank $r$, say) 
over $X$ such that  
\[
 h^0(\chi \otimes \Delta) = 0,
\]
i.e.,  $\chi \otimes \Delta$ has no nonzero global holomorphic sections. 
The condition implies that $\chi$ is necessarily semistable, and means 
that the (equivalence class of) $\chi$ lies on the complement of the 
generalized theta divisor in the moduli space of semistable vector 
bundles of rank $r$ and degree $0$ on $X$ (see \cite{newfay}, 
\cite{Seshadri} and \cite{Drezet}).  
It also follows immediately from Weil's criterion for flatness \cite{Gunning}
that $\chi$ is actually a flat vector bundle 
(see \cite{hip} page 275 for details). By definition, 
since $\chi$ is flat, sections $h$ of $\chi$ have the property that they 
lift to ${\bold C}^r$-vector functions $\widetilde{h}$ 
defined on the universal cover 
$\widetilde X$ of $X$ such that
\[
h(R \widetilde{p}) = \chi(R) h(\widetilde{p})
\]
for all $\widetilde{p} \in \widetilde{X}$ where $R$ is any element of the 
group of deck transformations Deck($\widetilde{X}/X)$ $ \cong \pi_1(X)$
and where $R \to \chi(R) \in GL(r)$ is a (constant) 
factor of automorphy associated with the bundle $\chi$.
(We somewhat abuse the notation denoting a constant factor of automorphy
and the corresponding flat vector bundle by the same letter.)
 
The main object of this section is  to define an object $K(\chi; \cdot, 
\cdot)$ associated with any such bundle $\chi$ which we shall 
call the {\it Cauchy kernel} for the bundle $\chi$. Let $M$ denote the 
Cartesian product $M=X \times X$ and let $\pi_1 \colon M \to X$ be the 
projection map onto the first coordinate and $\pi_2 \colon M \to X$ the 
projection onto the second coordinate.  The defining property of $K(\chi; 
\cdot, \cdot)$ is that $K(\chi; \cdot, \cdot)$ be a meromorphic 
mapping of the vector bundles $\pi_2^* \chi$ and $\pi_1^* \chi \otimes 
\pi_1^*\Delta \otimes \pi_2^* \Delta$ on $M$ which is holomorphic outside 
of the diagonal ${\cal D} = \{(p,p) \in M \colon p \in X\}$, where it has a 
simple pole with residue $I_r$.  More precisely, the latter condition 
means the following: for any $\widetilde{p}_0 \in \widetilde{X}$ and any local 
parameter $t$ on $X$, if we let $\sqrt{dt}$ be the corresponding local 
holomorphic frame for $\Delta$ lifted to a neighborhood of 
$\widetilde{p}_0$ on $\widetilde{X}$, then near $(\widetilde{p}_0, 
\widetilde{p}_0) \in \widetilde{X} \times \widetilde{X}$ the lift of 
$K(\chi; \cdot, \cdot)$ to $\widetilde{X} \times \widetilde{X}$ has the form
\[
\frac{ K(\chi; \widetilde{p},\widetilde{q}) }{ \sqrt{dt}(\widetilde{p}) 
\sqrt{dt}(\widetilde{q}) } = 
\frac{1}{ t(\widetilde{p}) - t(\widetilde{q}) } 
\left[ I_r + O\left(\sqrt{ |t(\widetilde{p})|^2 + |t(\widetilde{q})|^2 } 
\right) \right]. 
\]
Thus, if $e$ is in the fiber of $\chi$ at a point $q$ on $X$, and $t$ is 
a local parameter of $X$ centered at $q$, then $K(\chi; \cdot, q) 
\frac{e}{\sqrt{dt}(q)}$ is a meromorphic section of $\chi \otimes \Delta$ 
that has a single simple pole at $q$, with a residue (in terms of the 
local parameter $t$) equal to $e \sqrt{dt}(q)$.  Note that since $\chi 
\otimes \Delta$ has no nontrivial global holomorphic sections, such a 
meromorphic section is unique whenever it exists. 
Note that when $X$ is the Riemann sphere and $\chi$ is (necessarily) 
trivial, then
\[K(\chi;p,q) = \frac{I_r}{t(p) - t(q)} \sqrt{dt}(p) \sqrt{dt}(q)
\]
where $t$ is the standard coordinate on the complex plane, i.e., we get 
the usual Cauchy kernel.   
Existence was shown in \cite{hip} by exhibiting an explicit formula for 
$K(\chi; \cdot, \cdot)$; the construction involved  
using a representation of $X$ as the normalized Riemann surface for an 
algebraic curve $C$ embedded in ${\bold P}^2$ and representing the bundle $E 
= \chi \otimes \Delta \otimes {\cal O}(1)$ as the kernel bundle 
associated with a determinantal representation  of the curve $C$
\begin{gather} 
C = \{[\mu_0,\mu_1,\mu_2] \in {\bold P}^2 \colon \det (\mu_1\sigma_2 - \mu_2 
\sigma_1 + \mu_0 \gamma) = 0\} \notag \\
E(\mu) = \ker (\mu_1\sigma_2 - \mu_2 
\sigma_1 + \mu_0 \gamma).
\notag
\end{gather}

When the rank $r$ of the vector bundle $\chi$ is 1, one can get an 
explicit formula (in terms of the Abel-Jacobi map and classical theta 
functions on the Jacobian variety of $X$) for the Cauchy kernel (see 
\cite{hip}).  Details of this formula will be reviewed in Section 
\ref{S:linebundle} of this paper, where other special aspects of the line bundle 
case will also be discussed.

Our purpose in this section is to give an alternative proof of the existence 
of such a Cauchy kernel $K(\chi; \cdot, \cdot)$ for 
a flat vector bundle $\chi$ by a direct, simple, more abstract argument
(without relying on representing $X$ as the normalizing Riemann surface for a 
curve $C$ having a determinantal representation as in \cite{hip}).
This is the content of the following theorem.

\begin{theorem} \label{T:cauchyker}
Let $\chi$ be a flat vector bundle over the Riemann surface $X$ with 
$h^0(\chi \otimes \Delta) = 0$ as above.  Then the Cauchy kernel $K(\chi;
\cdot, \cdot)$ exists, i.e., there is a unique meromorphic 
mapping of the vector bundles $\pi_2^* \chi$ and $\pi_1^* \chi \otimes 
\pi_1^*\Delta \otimes \pi_2^* \Delta$ on $M=X \times X$ 
which is holomorphic outside 
of the diagonal ${\cal D} = \{(p,p) \in M \colon p \in X\}$, where it has a 
simple pole with residue $I_r$. 
\end{theorem}

\begin{pf}
Note that we have the following exact sequence of vector bundles over 
$M=X \times X$:  
\begin{equation} \label{exseq1}
0 \to {\cal O}(-{\cal D}) \to {\cal O} \to {\cal O}|_{\cal D} \to 0.
\end{equation}
For ease of notation, define vector bundles $F$ and $W$ over $M$ and $V$ 
and $K$ over $X$ by
\begin{gather}
F= \pi_1^* \chi \otimes \pi_1^* \Delta \otimes \pi_2^* \chi^\vee \otimes 
\pi_2^* \Delta \notag \\
W = F \otimes {\cal O}({\cal D}) \notag  \\
V = \chi \otimes \Delta \notag  \\
K = \text{ the canonical line bundle on } X \notag
\end{gather}
where $\chi^\vee$ is the dual bundle of $\chi$.
Tensoring the exact sequence \eqref{exseq1} with $W$ gives us
\begin{equation} \label{exseq2}
0 \to F \to W \to W \otimes {\cal O}|_{\cal D} \to 0.
\end{equation}
The map of taking the residue along the diagonal defines a linear mapping 
$$
{\cal R} \colon H^0(M,W) \rightarrow H^0(X, End\ V).
$$
Our goal is to 
show that there exists a unique element $K(\chi; \cdot, \cdot)$ of 
$H^0(M,W)$ so that ${\cal R}(K(\chi; \cdot, \cdot)) = I_V$.

Note that the bundle $W \otimes {\cal O}|_{\cal D}$ can be identified 
with the bundle $End\ V$ of endomorphisms of $V$. Moreover the residue 
mapping ${\cal R}$ is exactly the mapping from $H^0(M,W)$ into $H^0(X, 
End\ V)$ induced by the mapping $W \rightarrow \left. W \otimes
{\cal O}\right|_{\cal D}$ in \eqref{exseq2}. 

The vector bundle exact sequence \eqref{exseq2} 
induces  (see the Basic Fact 
on page 40 of \cite{gh}) the exact cohomology sequence
\begin{align} \notag
0 \to H^0(M,F) & \to H^0(M,W) \to H^0(X, End\ V) \to \\
   \to H^1(M, F) & \to \cdots. 
   \label{exseq3}
\end{align}

Next we argue that (i) $h^0(F) = 0$ and (ii) $h^1(F) = 0$.  The 
statement (i) follows easily from our assumption that $h^0(V) = 0$. 
 As for statement (ii) it follows from the Kunneth 
formulas (see page 58 of \cite{gh}) that
\begin{align}
H^1(M,F) = & H^1(M, \pi_1^*(V) \otimes \pi_2^*( V^\vee 
\otimes K)) 
\notag \\
\cong & \left(H^0(M, \pi_1^*(V)) \otimes H^1(M, \pi_2^*(V^\vee 
\otimes K)) \right)  \\
 & \oplus \left( H^1(M, \pi_1^*(V)) \otimes 
H^0(M, \pi_2^*(V^\vee \otimes K) ) \right).  \label{kunneth}
\end{align}
By our assumption that $h^0(V)=0$ it follows that the 
first term on the right hand side of \eqref{kunneth} is 0.  Since 
$h^0(V) = 0 $ we also have $h^0(V^\vee \otimes 
K) = 0$ as well, by the Riemann-Roch Theorem for vector bundles on an 
algebraic curve (see \cite{Gunning}) and the assumption that deg $V=r(g-1)$.
Hence the second term in \eqref{kunneth} is zero as 
well.  This verifies the desired fact (ii).

Hence the exact sequence \eqref{exseq3} collapses to
\begin{equation} \label{exseq4}
0 \to H^0(M,W) \overset{{\cal R}}{\to} H^0(X,End\ V) \to 0.
\end{equation}
It follows that ${\cal R} \colon H^0(M,W) \to H^0(X,End\ V)$
is an isomorphism and the Theorem follows.
\end{pf}

Let $\widetilde{p}_0$ be an arbitrary point of $\widetilde{X}$, and let 
$t$ be a local coordinate for $\widetilde{X}$ near $\widetilde{p}_0$.  
Then by definition the Cauchy kernel $K(\chi; \cdot, \cdot)$ is such that
\[ 
\left(t(\widetilde p) - t(\widetilde q) \right)
\frac{K(\chi; \widetilde p, \widetilde q)}{\sqrt{dt}(\widetilde p) 
\sqrt{dt}(\widetilde q)}
\]
is analytic in $(\widetilde p, \widetilde q)$ near $(\widetilde p_0, 
\widetilde p_0)$ with value at $(\widetilde p_0, \widetilde p_0)$ equal 
to $I_r$; hence, for $(\widetilde p, \widetilde q)$ close to $(\widetilde 
p_0, \widetilde p_0)$, $K(\chi; \cdot, \cdot)$ has a representation of 
the form
\begin{multline} \label{cauchyexp}
\frac{K(\chi; \widetilde p, \widetilde q)}{\sqrt{dt}(\widetilde p) 
\sqrt{dt}(\widetilde q)} \\
 = \frac{1}{t(\widetilde p) - t(\widetilde q)}
\left[ I_r + \frac{A_\ell}{dt}(\widetilde p_0) t(\widetilde p) +
\frac{A}{dt}(\widetilde p_0) t(\widetilde q) + O\left( |t(\widetilde 
p)|^2 + |t(\widetilde q)|^2 \right) \right]
\end{multline}
for appropriate holomorphic matrix $K$-valued coefficients $A(\widetilde 
p_0)$ and $A_\ell(\widetilde p_0)$.  An important point for us is the 
following result.

\begin{lemma} \label{L:cauchyexp}
Let $A_\ell$ and $A$ be the linear coefficients appearing 
in the Laurent expansion of the Cauchy kernel $K(\chi; \cdot, \cdot)$ as 
in \eqref{cauchyexp}.  Then
\[A(\widetilde p_0) + A_\ell(\widetilde p_0) = 0
\]
for all $\widetilde p_0 \in \widetilde X$.
\end{lemma}

\begin{pf}
Let $\widetilde{p}_0$ be an arbitrary point of $\widetilde X$ and let $t$ 
be a local coordinate for $\widetilde X$ near $\widetilde p_0$.  For 
$\widetilde p \in X$ near $\widetilde p_0$, define
$$
f(\widetilde p) = (t(\widetilde p) - t(\widetilde q)) \cdot
\left. \frac{K(\chi; \widetilde p, \widetilde q)}{\sqrt{dt}(\widetilde p) 
\sqrt{dt}(\widetilde q)} \right|_{\widetilde p = \widetilde q}.
$$
From \eqref{cauchyexp} we see that
$$
f(\widetilde p) = I_r + \left[ \frac{A_\ell}{dt}(\widetilde p_0) + 
\frac{A}{dt}(\widetilde p_0)\right] t(\widetilde p) + O(|t(\widetilde p)|^2).
$$
In particular,
\begin{equation} \label{der1}
\left. \frac{df}{dt}(\widetilde p) \right|_{\widetilde p = \widetilde p_0} =
\frac{A_\ell}{dt}(\widetilde p_0) + \frac{A}{dt}(\widetilde p_0).
\end{equation}
On the other hand, we can use $t'(\widetilde p') = t(\widetilde p') - 
t(\widetilde p)$ as a local coordinate for the variable $\widetilde p'$
near the point $\widetilde p \in X$.  From \eqref{cauchyexp} again we have
\begin{multline} 
[(t(\widetilde p') - t(\widetilde p)) - (t(\widetilde q') - t(\widetilde p))]
\frac{K(\chi; \widetilde p', \widetilde q')}{\sqrt{dt}(\widetilde p')
\sqrt{dt}(\widetilde q')}  \\
\notag = I_r + \frac{A_\ell}{dt}(\widetilde p) (t(\widetilde p') - 
t(\widetilde p)) + \frac{A}{dt}(\widetilde p) (t(\widetilde q') - 
t(\widetilde p)) + O(|t(\widetilde p') - t(\widetilde p)|^2 + 
|t(\widetilde q') - t(\widetilde p)|^2).  
\end{multline}
Evaluation of both sides of this equation at $\widetilde p' = \widetilde 
q' = \widetilde p$ yields $f(\widetilde p) = I_r$ from which we get
\begin{equation} \label{der2}
\frac{df}{dt}(\widetilde p) = 0
\end{equation}
for all $\widetilde p \in X$.  Comparison of \eqref{der1} and 
\eqref{der2} now gives $A(\widetilde p_0) + A_\ell(\widetilde p_0) = 0$ 
as asserted.
\end{pf}

We can use these coefficients $A(\widetilde p)$ and $A_\ell(\widetilde q)$ 
defined by \eqref{cauchyexp} to define connections $\nabla_\chi$ on 
$\chi$ and $\nabla^*_\chi$ on $\chi^\vee$ according to the formulas
\begin{align}
\nabla_\chi y &= A y + dy  \notag \\
\nabla^*_\chi x &= A_\ell^T x + dx. \label{connection}
\end{align}
for local holomorphic sections $y$ of $\chi$ and $x$ of $\chi^\vee$.
The result $A_\ell + A = 0$ from Lemma \ref{L:cauchyexp} is equivalent to the 
fact that $\nabla_\chi$ and $\nabla_\chi^*$ are {\it dual to each other}, i.e.
\[
d(x^Ty) = x^T (\nabla_\chi y) + (\nabla_\chi^* x)^T y
\]
for local holomorphic sections $y$ of $\chi$ and $x$ of $\chi^\vee$, 
where $(y, x) \to x^Ty$ is the pairing between $\chi$ and $\chi^\vee$.
Moreover, from the formula for $\nabla_\chi$ we see that the connection 
matrix associated with $\nabla_\chi$ is a matrix of holomorphic 
$(1,0)$-forms. 
Hence, $\nabla_\chi$ is {\it compatible with the complex structure of}
$X$ and moreover, since we are in complex dimension 1, the connection 
$\nabla_\chi$ is {\it flat}, i.e., $\nabla_\chi$ has {\it zero curvature} 
(see Section~5 of Chapter~0 of \cite{gh} for all relevant definitions).  
The existence of such a flat connection on $\chi$ in turn implies that 
$\chi$ itself is a flat vector bundle (see \cite{mst} pages 294--295).
In general there are many choices of distinct flat connections on a flat 
vector bundle; our construction via the Cauchy kernel provides a canonical 
choice of such a flat connection (up to a choice of a bundle $\Delta$ 
of half-order 
differentials).  An explicit formula for $\nabla$ in the line 
bundle case is given in Section \ref{S:linebundle}.

{\bf Remark:} The proof of Theorem \ref{T:cauchyker}
used only the fact that $\deg\chi = 0$ and did not use the flatness of $\chi$.
Since  the existence of a flat connection (compatible with the 
complex structure) implies that the 
bundle is flat, our construction 
gives a direct proof of 
the flatness of $\chi$ independent of Weil's theorem.

\section{The abstract interpolation problem} \label{S:absint}

In this section we consider as given two flat vector bundles $\chi$ and 
$\widetilde \chi$ over the Riemann surface $X$ for which both $h^0(\chi 
\otimes \Delta)=0$ and $h^0(\widetilde \chi \otimes \Delta)=0$, where 
again, $\Delta$ is a line bundle of half-order differentials over $X$.  
We are interested in studying pole-zero interpolation conditions imposed 
on a bundle map of $\chi$ to $\widetilde \chi$.  The data for the 
interpolation problem is as follows.  We assume that we are given 
$n_\infty$ distinct points $\mu^1, \dots, \mu^{n_\infty}$ (the 
prescribed poles) together with $n_0$ distinct points $\lambda^1, \dots, 
\lambda^{n_0}$ (the prescribed zeros). For each fixed index $j$ ($j=1, 
\dots, n_\infty$) we specify a linearly independent set $\{u_{j1}, \dots, 
u_{j,s_j}\}$ of $s_j$ vectors in the fiber $\widetilde \chi(\mu^j)$ of 
$\widetilde \chi$ over $\mu^j$ (the prescribed pole vectors) and for each 
fixed index $i$ ($i=1, \dots, n_0$) we specify a linearly independent set 
$\{x_{i1}, \dots, x_{i,t_i}\}$ of $t_i$ vectors in the fiber $\widetilde 
\chi^\vee(\lambda^i)$ of the dual bundle $\widetilde \chi^\vee$ of 
$\widetilde \chi$ (the prescribed null vectors).  Also, for each pair of 
indices ($i,j$) for which $\lambda^i = \mu^j=:\xi^{ij}$, we specify a 
collection $\{\rho_{ij,\alpha\beta} \colon 1 \le \alpha \le t_i, 1 \le 
\beta \le s_j\}$ of numbers that depend on the choice of the local 
parameter at the point $\xi^{ij}$.  The {\it Abstract Interpolation 
Problem} (ABSINT) which we study in this section is the following: 
{\it determine if 
there exists a bundle map $T \colon \chi \to \widetilde \chi$ with 
transpose $T^\vee \colon \widetilde \chi^\vee \to \chi^\vee$ such that:
\begin{enumerate}
\item[(i)] $T$ has poles only at the points $\{\mu^1, \dots, \mu^{n_\infty}\}$; 
for each $j=1, \dots, n_\infty$, the pole of $T$ at $\mu^j$ is simple, 
and the residue $R_j = \text{ Res }_{p=\mu^j} T \colon \chi(\mu^j) \to 
\widetilde \chi(\mu^j)$ of $T$ at $\mu^j$ is such that $\{u_{j1}, \dots, 
u_{j,s_j}\}$ spans the image space $\text{im } R_j$ of $R_j$.
\item[(ii)]  The bundle map $(T^\vee)^{-1} \colon \chi^\vee \to \widetilde 
\chi^\vee$ has poles only at $\{\lambda^1, \dots, \lambda^{n_0}\}$; for 
each $i=1, \dots, n_0$, the pole of $(T^\vee)^{-1}$ at $\lambda^i$ is 
simple and the residue $\widehat R_i = \text{ 
Res}_{p=\lambda^i}(T^\vee)^{-1} \colon \chi^\vee(\lambda^i) \to 
\widetilde \chi^\vee (\lambda^i)$ of $(T^\vee)^{-1}$ at $\lambda^i$ is 
such that $\{x_{i1}, \dots, x_{i,t_i}\}$ spans the image space $\text{ im 
} \widehat R_i$ of $\widehat R_i$.
\item[(iii)] For each pair of indices ($i,j$) for which $\lambda^i = \mu^j =: 
\xi^{ij}$, and for $\alpha = 1, \dots, t_i$, let $x_{i \alpha}(p)$ be a 
local holomorphic section of $\widetilde \chi^\vee$ with $x_{i 
\alpha}(\xi^{ij}) = x_{i \alpha}$ such that $T^\vee(p)x_{i \alpha}(p)$ 
has analytic continuation to $p = \xi^{ij}$ with value at $p = \xi^{ij}$ 
equal to 0.  Then
\[
(\nabla^*_{\widetilde \chi} x_{i \alpha}(\xi^{ij}))^T u_{j \beta} = \rho_{ij, 
\alpha \beta}
\]
for $\beta = 1, \dots , s_j$.
\end{enumerate}
 When such a bundle map $T$ exists, give an 
explicit formula for the construction of $T$.}

In order for solutions to exist, the compatibility condition
\begin{equation}  \label{comp}
x_{i \alpha}u_{j \beta} = 0
\text{ whenever } \lambda^i = \mu^j.
\end{equation}
must hold.  This follows from the requirement that the meromorphic 
section 
$$
x_{i \alpha}(p)T(p)
$$ 
be analytic at the point $p=\xi^{ij}:= 
\lambda^i = \mu^j$.  Hence we shall always assume that our data 
collection 
\begin{equation}  \label{dataset}
{\boldsymbol \omega} = \{(x_{i\alpha},\lambda^i), (u_{j\beta}, \mu^j), 
\rho_{ij, \alpha \beta} \}
\end{equation}
also satisfies this compatibility condition.

It will be convenient to work with an alternate form of the interpolation 
condition (iii) in (ABSINT).  Suppose that $u(p) = T(p) \varphi(p)$, 
where $\varphi$ is a local holomorphic section of $\chi$ near $\xi^{ij}$ 
chosen so that $\text{Res}_{p=\xi^{ij}} u(p) = u_{j \beta}$ with respect 
to the local coordinate $t^{ij}$ centered at $\xi^{ij}$.  Then
\begin{align}
 \left(x_{i \alpha}(p) \right)^T\left( t^{ij}(p) u(p) \right)
& = t^{ij}(p) \cdot \left(x_{i \alpha}(p)\right)^T T(p) \varphi(p) \\
& = \left( t^{ij}(p)  T^\vee(p) x_{i \alpha}(p)\right)^T \varphi(p)
\end{align}
has a double order zero at $p=\xi^{ij}$, and hence
\[
d\left( x_{i \alpha}(p)^T (t^{ij}(p) u(p) ) \right)|_{p=\xi^{ij}} = 0.
\]
Since $\nabla_{\widetilde{\chi}}^*$ and $\nabla_{\widetilde{\chi}}$ are 
dual connections as a consequence of Lemma \ref{L:cauchyexp}, we 
therefore have
\[
\left( \nabla_{\widetilde{\chi}}^* x_{i \alpha}(\xi^{ij})\right)^T u_{j 
\beta} + x_{i \alpha}^T \nabla_{\widetilde{\chi}} (t^{ij}(p) u(p))|_{ p= 
\xi^{ij}} = 0.
\]
Thus the interpolation condition in part (iii) of (ABSINT) can be 
expressed alternatively as
\begin{equation} \label{coupledint}
x_{i \alpha}^T \nabla_{\widetilde{\chi}}(t^{ij}(p) u(p))|_{p = \xi^{ij}} = 
- \rho_{ij, \alpha \beta}
\end{equation}
where $u(p) = T(p) \varphi(p)$ for a local holomorphic section 
$\varphi$ of $\chi$ 
near $\xi^{ij}$ such that $\text{Res}_{p=\xi^{ij}}u(p) = u_{j \beta}$.

In the scalar case ($r=1$), the compatibility condition \eqref{comp} can 
never be satisfied in a nontrivial way and hence the third set of 
interpolation conditions is absent under our assumptions; this 
corresponds to the fact that a scalar meromorphic function cannot have a 
zero and a pole at the same point $\xi^{ij}$.  Moreover, for the case 
$r=1$, necessarily $t_i=1$ for all $i$ and $s_j=1$ for all $j$.  In the 
case where both $\chi$ and $\widetilde{\chi}$ are trivial (or more 
generally if we use coordinates with respect to a local holomorphic frame 
for $\widetilde{\chi}$ near $\mu^j$ or $\widetilde{\chi}^\vee$ near 
$\lambda^i$), there is no loss of generality in taking $x_i:=x_{i1}=1$ 
and $u_j:=u_{j1}=1$ for all $i$ and $j$.  Thus the only remaining 
relevant data are the zeros $\lambda^1, \dots \lambda^{n_0}$ 
and the poles $\mu^1, \dots, \mu^{n_\infty}$ (all assumed here to 
be distinct).  As is standard in algebraic geometry, the formal sum
\[
{\boldsymbol \lambda} - {\boldsymbol \mu}:= \lambda^1 + \dots + 
\lambda^{n_0} - \mu^1 
- \dots - \mu^{n_\infty}
\]
is said to be a {\it divisor} on $X$.   If $f$ is a meromorphic function, 
the associated principal divisor $(f)$ is defined to be the formal sum
$ p^1 + \dots + 
 p^{n_0} - q^1 - \dots - q^{n_\infty}$ where the $p^i$'s are the zeros of 
$f$ and the $q^j$'s are the poles of $f$ (with repetitions according to 
respective multiplicities).  Associated with any divisor ${\boldsymbol 
\lambda} - {\boldsymbol \mu}$ as above is the vector bundle 
${\cal O}({\boldsymbol \mu} - 
{\boldsymbol \lambda})$ 
whose holomorphic sections can be identified 
with global 
meromorphic functions $h$ such that
\[
(h) \ge {\boldsymbol \lambda} - {\boldsymbol \mu},
\] 
i.e., such that the zeros of $h$ include the points $\lambda^1, \dots, 
\lambda^{n_0}$ (all with multiplicity at least 1) and the poles of $h$ 
are a subset of $\mu^1, \dots, \mu^{n_\infty}$ (all with multiplicity at 
most 1).

It is convenient for us to introduce matrix analogues of these ideas.  
Let $\boldsymbol \omega$ be an interpolation data set as in \eqref{dataset}.
Let us introduce the notation
\begin{equation} \label{poledata}
({\boldsymbol \mu}, {\bold u}) = \{ (\mu^j, u_{j \beta}) \colon 1 \le j \le 
n_\infty, 1 \le \beta \le s_j \}
\end{equation}
for the pole part of $\boldsymbol \omega$. 
We let ${\cal M}(\widetilde{\chi}\otimes \Delta)$ be the sheaf of 
meromorphic sections of $\widetilde{\chi} \otimes \Delta$.
 For $U$ an open subset of $X$ we 
define
\begin{align} \notag
{\cal O}(\widetilde{\chi} \otimes \Delta)({\boldsymbol \mu}, {\bold u}) (U) = &
\{u \in {\cal M}(\widetilde{\chi} \otimes \Delta ) (U) \colon u 
\text{ has poles only at } \mu^j, \\
& u_{-1}:=\text{Res}_{p=\mu^j} \in \text{span }\{u_{j \beta} \colon 1 \le 
\beta \le s_j\} \}
\notag
\end{align}
and
\begin{align} \notag
{\cal O}(\widetilde{\chi}\otimes \Delta)({\boldsymbol \omega})(U) = 
& \{u \in {\cal O}(\widetilde{\chi} \otimes \Delta) ({\boldsymbol \mu}, {\bold 
u})(U) \colon \text{ if } \lambda^i \ne \mu^j, \text{ then } x_{i 
\alpha}u(\lambda^i) = 0; \\
\notag
& \text{if } \lambda^i = \mu^j =:\xi^{ij}, \text{ there is a local 
holomorphic section } x_{i \alpha}(p) \text{ of } \widetilde{\chi}^\vee \\
\notag
& \text{such that (i) } x_{i \alpha}(\xi^{ij}) = x_{i \alpha}, \text{ (ii) }
x_{i \alpha}(p)^T \frac{u}{ \sqrt{ dt^{ij} } }(p) \text{ has analytic} \\
\notag
& \text{continuation to } p = \xi^{ij} 
 \text{with value } 0 \text{ there, and} \\
 \notag
 & \text{ (iii) }
x_{i \alpha}^T \nabla_{\chi}\left(t^{ij} \frac{u}{ \sqrt{dt^{ij}} }\right)
= - \sum_{\beta = 1}^{s_j}\rho_{ij,\alpha \beta}
 c_\beta \frac{1}{ \sqrt{dt^{ij}} (\xi^{ij})} \\
\notag 
& \text{if Res}_{p=\mu^j} u = \sum_{\beta=1}^{s_j} u_{j \beta} c_\beta.\}
\end{align}
It is obvious that ${\cal O}(\widetilde{\chi} \otimes 
\Delta)({\boldsymbol \mu}, {\bold u})$ and ${\cal O}(\widetilde{\chi} 
\otimes \Delta)({\boldsymbol \omega})$ are locally free sheaves of rank 
$r$, and we denote the corresponding rank $r$ vector bundles by 
$(\widetilde{\chi} \otimes \Delta)({\boldsymbol \mu}, {\bold u})$ and 
$(\widetilde{\chi} \otimes \Delta)({\boldsymbol \omega})$.  It is also 
obvious that $T$ is a solution of (ABSINT) if and only if $T$ is an 
isomorphism from $\chi \otimes \Delta$ to $(\widetilde{\chi} \otimes 
\Delta)({\boldsymbol \omega})$.

The solution of the zero-pole interpolation problem introduced at the 
beginning of this section is as follows.

\begin{theorem} \label{T:absint}
Define a $n_0 \times n_\infty$ block matrix $\Gamma = [\Gamma_{ij}]$
($1 \le i \le n_0$, $1 \le j \le n_\infty$) 
where the block entry $\Gamma_{ij}$ in turn is a $t_i \times s_j$ matrix
$\Gamma_{ij} = [\Gamma_{ij, \alpha \beta}]$ ($1 \le \alpha \le t_i$, $1 
\le \beta \le s_j$) with matrix entries $\Gamma_{ij, \alpha \beta}$ given by
\begin{equation} \label{defGamma}
\Gamma_{ij, \alpha \beta} = 
\begin{cases} - x_{i \alpha}^T K(\widetilde{\chi};\lambda^i, \mu^j) u_{j 
\beta}, &\text{if } \lambda^i \ne \mu^j; \\
-\rho_{ij, \alpha \beta} & \text{if } \lambda^i = \mu^j.
\end{cases}
\end{equation}
In addition we introduce the block matrices
\begin{gather} \notag
{\bold u}_i = \begin{bmatrix} u_{i1} & \dots & u_{i s_i} \end{bmatrix}, 
\quad
{\bold x}_j = \begin{bmatrix} x_{j1} & \dots & x_{j t_i} \end{bmatrix}, \\
K_{{\boldsymbol \mu}, {\bold u}}(p) = 
\begin{bmatrix} K(\widetilde{\chi}; p, \mu^1) {\bold u}_1 & \dots & 
K(\widetilde{\chi}; p, \mu^{n_\infty}) {\bold u}_{n_\infty} \end{bmatrix}, \\
K^{{\bold x},{\boldsymbol \lambda}}(q) = \begin{bmatrix} {\bold x}^T_1 
K(\widetilde{\chi}; \lambda^1, q) \\
\vdots \\ {\bold x}^T_{n_0} K(\widetilde{\chi}; \lambda^{n_0}, q)\end{bmatrix}.
\end{gather}
Let $q$ be a point of $X$ disjoint from all the interpolation nodes 
$\lambda^1, \dots, \lambda^{n_0}$, $\mu^1, \dots, \mu^{n_\infty}$
and let $Q$ be an invertible linear map of the fiber 
space $\chi(q)$ to the fiber 
space $\widetilde{\chi}(q)$.
Then the abstract interpolation problem (ABSINT) has a solution $T$ with 
value $Q$ at the point $q$ if and only if the matrix $\Gamma$ is square 
and invertible and
\begin{equation} \label{residues}
[K(\widetilde{\chi}; p^i,q) + K_{{\boldsymbol\mu}, {\bold u}}(p^i) \Gamma^{-1}
K^{{\bold x},{\boldsymbol \lambda}}(q)] 
Q (\text{Res}_{p^i}\ K(\chi; \cdot , q)^{-1}) = 0
\end{equation}
at each pole $p^i$ of $K(\chi; \cdot, q)^{-1}$.
In this case the unique solution $T$ of the interpolation problem (ABSINT)
with value 
$Q$ at $q$ is given by
\begin{equation}  \label{solution}
T(p) = [K(\widetilde{\chi}; p,q) + K_{{\boldsymbol \mu}, {\bold u}}(p) 
\Gamma^{-1} K^{{\bold x},{\boldsymbol \lambda}}(q)] Q K(\chi; p,q)^{-1}
\end{equation}
with inverse given by
\begin{equation} \label{inversesolution}
T^{-1}(p) = K(\chi;p,q)^{-1} T^{-1}(q)[K(\widetilde{\chi};q,p) + 
K_{{\boldsymbol \mu},{\bold u}}(q) \Gamma^{-1} K^{ {\bold x},
{\boldsymbol \lambda} }(p)].
\end{equation}
\end{theorem}

Two special cases of formula  \eqref{solution} deserve to be mentioned.  
The first is the  case where $n_0=n_\infty = 1$, $t_1=s_1=1$ and 
$\lambda^1 \ne \mu^1$.  If we set $x =x_1$, $\lambda = \lambda^1$, $\mu 
= \mu^1$ and $u = u_1$, then $\Gamma = -x K(\widetilde{\chi}; \lambda, 
\mu) u$ is just a number and the formula \eqref{solution} becomes
 \begin{equation} \label{3.12a}
T(p) K(\chi; p,q) T(q)^{-1} = 
\dfrac{ K(\widetilde{\chi}; p,q) - K(\widetilde{\chi}; p,\mu) u x^T 
K(\widetilde{\chi}; \lambda,q)}{x^T K(\widetilde{\chi}; \lambda, \mu) u}.
\end{equation}
In the line bundle case, the identity \eqref{3.12a} reduces to the Fay trisecant 
identity and will be discussed in Section \ref{S:linebundle}.  The second
special case of interest 
is the case where the given zero and pole vectors at each 
interpolation node span the whole fiber space.  In this case there is an 
explicit multiplicative formula for the interpolant $T$ 
in terms of the prime form 
$E(p,q)$; this will be discussed in detail at the end of Section 
\ref{S:linebundle}.

A second version of the abstract interpolation problem (ABSINT) has the 
same form (i), (ii) and (ii) as (ABSINT), but with the input bundle 
$\chi$ left also as an unknown to be found, subject to the proviso that 
it also be flat and have $h^0(\chi \otimes \Delta)=0$.  This version of 
the problem was studied in \cite{hip} in a more concrete setting where 
$X$ is the normalizing Riemann surface for an algebraic curve $C$ 
embedded in ${\bold P}^2$ having a maximal rank $r$
 determinantal representation; we will discuss the connections of this 
 setup with ours in Sections 5 and 6.  At this time we also state the 
 solution to the modified (ABSINT).
 
\begin{theorem} \label{T:bvabsint}
Let $\boldsymbol \omega$ be a data set for (ABSINT) 
as above and form the matrix 
$\Gamma$ as in \eqref{defGamma}.  Then there exists a flat bundle 
$\chi$ with $h^0(\chi \otimes \Delta) = 0$ and a meromorphic bundle map 
$T \colon \chi \to \widetilde{\chi}$ satisfying the interpolation 
conditions (i), (ii) and (iii) of (ABSINT) if and only if $\Gamma$ is 
square and invertible.
\end{theorem}

To prove Theorem \ref{T:absint} we need some preliminary lemmas.

\begin{lemma}  \label{L1:absint}
A global meromorphic section of $\widetilde{\chi}\otimes \Delta$ is in 
${\cal O}(\widetilde{\chi} \otimes \Delta) ({\boldsymbol \mu},{\bold u})(X)$ 
if and only if $h$ has the form
\[
h= \sum_{j=1}^{n_\infty} \sum_{\beta = 1}^{s_j} K(\widetilde{\chi}; p, 
\mu^j) u_{j \beta} c_{j \beta}
\]
for some scalars $c_{j \beta}$.
\end{lemma}

\begin{pf}  Suppose $h \in {\cal O}(\widetilde{\chi} \otimes \Delta) 
( {\boldsymbol \mu}, {\bold u})(X)$.  Choose scalars $c_{j \beta}$ so that
\[
\text{Res}_{p = \mu^j} h(p) = \sum_{\beta=1}^{s_j} u_{j \beta} c_{j \beta}
\]
and set
\[
\widehat{h}(p) = \sum_{j=1}^{n_\infty} \sum_{\beta=1}^{s_j} 
K(\widetilde{\chi}; p, \mu^j) u_{j \beta} c_{j \beta}.
\]
Then $\widehat{h} \in {\cal M}(\widetilde{\chi} \otimes \Delta)(X)$ and $h 
- \widehat{h} \in {\cal O}(\widetilde{\chi} \otimes \Delta)(X)$. Thus 
$h = \widehat{h}$ since ${\cal O}(\widetilde{\chi} \otimes \Delta)(X) 
= H^0(X, \widetilde{\chi} \otimes \Delta) = 0$ by our standing 
assumptions on $\widetilde{\chi}$.  The converse direction follows easily 
from the defining properties of the Cauchy kernel $K(\widetilde{\chi}; 
\cdot, \cdot)$.
\end{pf}

\begin{lemma} \label{L2:absint}
The map
\[
[[c_{j \beta}]_{1 \le \beta \le s_j}]_{1 \le j \le n_\infty} \to
h(p)=\sum_{j=1}^{n_\infty} \sum_{\beta = 1}^{t_j} K(\widetilde{\chi}; p, 
\mu^j) u_{j \beta} c_{j \beta}
\]
establishes a one-to-one correspondence between $\ker \Gamma$ and 
${\cal O}(\widetilde{\chi} \otimes \Delta)( {\boldsymbol \omega})(X)$.  
In particular
\[
\dim \ker \Gamma = h^0( (\widetilde{\chi} \otimes \Delta) 
({\boldsymbol \omega}) ).
\]
\end{lemma}

\begin{pf}
Suppose first that $h \in {\cal O}(\widetilde{\chi} \otimes \Delta)( 
{\boldsymbol \omega})(X)$.  In particular 
$h \in {\cal O}(\widetilde{\chi} \otimes 
\Delta) ({\boldsymbol \mu}, {\bold u})(X)$, so by Lemma \ref{L1:absint} there 
exists a collection of complex numbers $\{c_{j \beta}\}_{1\le j \le 
n_\infty, 1 \le \beta \le s_j}$ so that 
\[
h(p) = \sum_{j=1}^{n_\infty} \sum_{\beta=1}^{t_j} K(\widetilde{\chi}; p, 
\mu^j) u_{j \beta} c_{j \beta}.
\]
We next see what conditions the other requirements on $h$ for admission 
to the class ${\cal O}(\widetilde{\chi} \otimes \Delta) ({\boldsymbol 
\omega})(X)$ impose on the scalars $c_{j \beta}$.  

If $i$ is any index for which $\lambda^i \ne \mu^j$ for all $j$, we must have
\begin{align} \notag
0 = x_{i \alpha}h(\lambda^i) & = \sum_{j=1}^{n_\infty} \sum_{\beta = 1} 
^{t_j} x_{i \alpha} K(\widetilde{\chi}; \lambda^i, \mu^j) u_{j \beta} 
c_{j \beta} \\
\label{kerGamma1}
&= - \sum_{j=1}^{n_\infty} \sum_{\beta = 1}^{t_j} \Gamma_{ij, \alpha, 
\beta} c_{j \beta}.
\end{align}
If, on the other hand, $i$ is an index such that $\lambda^i = \mu^j =: 
\xi^{ij}$ for some index $j$, then we write the Laurent series expansion 
for $\frac{h}{\sqrt{dt^{ij}}}$  near $\xi^{ij}$ (where $t^{ij}(p)$ is a 
local coordinate for $X$ centered at $\xi^{ij}$)
\[
\frac{h}{ \sqrt{dt^{ij}} }= \left[ \frac{h}{\sqrt{dt^{ij}}}\right]_{-1} 
\frac{1}{t^{ij}} + \left[\frac{h}{\sqrt{dt^{ij}}}\right]_0 + O(|t^{ij}(p)|).
\]
We compute
\begin{equation} \label{res}
\left[\frac{h}{\sqrt{dt^{ij}}}\right]_{-1} = \text{Res}_{p=\xi^{ij}}\
\frac{h}{\sqrt{t^{ij}}}(p)
=\sum_{\beta=1}^{s_j} u_{j \beta} \sqrt{dt^{ij}}(\xi^{ij}) c_{j \beta}
\end{equation}
and hence
\begin{equation} \label{rescoeff}
\text{Res}_{p=\xi^{ij}} \frac{h}{\sqrt{dt^{ij}}}(p) = \sum_{\beta = 
1}^{s_j} u_{j \beta} {\bold c}_{j \beta}
\end{equation}
where ${\bold c}_{j \beta} = \sqrt{dt^{ij}}(\xi^{ij}) c_{j \beta}$. 
Moreover,, we have
\begin{align} \notag
\left[ \frac{h}{\sqrt{dt^{ij}}}\right]_0 = & \sum_{k \ne j} \sum_{\beta = 
1} ^{s_k} \frac{K(\widetilde{\chi}; \xi^{ij}, 
\mu^k)}{\sqrt{dt^{ij}}(\xi^{ij})} u_{k \beta}c_{k \beta} \\
\label{0coeff}
&+ \sum_{\beta=1}^{s_j} \frac{A_{\ell}(\xi^{ij})}{dt(\xi^{ij})} 
\sqrt{dt^{ij}}(\xi^{ij}) u_{j \beta} c_{j \beta}.
\end{align}
We next compute
\begin{gather} \notag
x_{i \alpha}^T ( 
A(\xi^{ij})\left[\frac{h}{\sqrt{dt^{ij}}}\right]_{-1} +  
\left[ \frac{h}{ \sqrt{dt^{ij}} } \right]_0 dt^{ij}(\xi^{ij}) )  =
\sum_{\beta=1}^{s_j} x^T_{i \alpha} A(\xi^{ij}) u_{j \beta} 
\sqrt{dt^{ij}}(\xi^{ij}) c_{j \beta} \\
\notag
 +\sum_{k \ne j} \sum_{\beta=1}^{s_j} x^T_{i \alpha} 
\frac{ K(\widetilde{\chi}; \xi^{ij}, \mu^k) }{ \sqrt{dt^{ij}}(\xi^{ij}) } 
u_{j \beta} c_{j \beta} dt^{ij}(\xi^{ij}) 
 + \sum_{\beta=1}^{s_j} x_{i \alpha} 
\frac{ A_\ell(\xi^{ij}) }{ \sqrt{dt^{ij}}(\xi^{ij}) } u_{j \beta} c_{j \beta} 
dt^{ij}(\xi^{ij}) \\
\notag
=  \left\{ \sum_{\beta=1}^{s_j} x_{i \alpha}^T (A(\xi^{ij}) + 
A_\ell(\xi^{ij})) u_{j \beta} c_{j \beta} + \sum_{k \ne j} 
\sum_{\beta=1}^{s_k} x_{i \alpha}^T K(\widetilde{\chi}; \xi^{ij}, \mu^k) 
u_{j \beta} c_{j \beta} \right\} \sqrt{dt^{ij}}(\xi^{ij}).
\end{gather}
By Lemma \ref{L:cauchyexp} the first term in the braces vanishes.  
Combining this fact with the formula \eqref{rescoeff} for the 
coefficients ${\bold c}_{j \beta}$, we see that the interpolation condition
\[
x_{i \alpha}^T \nabla_{\widetilde{\chi}}\left( t^{ij} 
\frac{h}{\sqrt{dt^{ij}}}\right) = - \sum_{\beta=1}^{s_j} \rho_{ij, \alpha 
\beta} {\bold c}_{j \beta}
\]
becomes
\[
\left\{ \sum_{k \ne j} \sum_{\beta = 1}^{s_k} x_{i \alpha}^T 
K(\widetilde{\chi}; \xi^{ij}, \mu^k) u_{j \beta} c_{j \beta} \right\}
\sqrt{dt^{ij}}(\xi^{ij}) =
- \left( \sum_{\beta = 1} ^{s_j} \rho_{ij, \alpha \beta} c_{j \beta} \right)
\sqrt{dt^{ij}}(\xi^{ij}).
\]
After canceling off $\sqrt{dt^{ij}}(\xi^{ij})$ and recalling that 
$\xi^{ij} = \lambda^i$ we see that this can be rewritten as
\begin{equation} \label{kerGamma2}
\sum_{j=1}^{n_\infty} \sum_{\beta=1}^{s_j} \Gamma_{ij, \alpha, \beta} c_{j 
\beta} = 0
\end{equation}
and this equation holds for all pairs of indices $(i, \alpha)$ such that 
$\lambda^i = \mu^j$ for some $j$.  Combining \eqref{kerGamma2} and 
\eqref{kerGamma1} we see that the column vector $[[c_{j \beta}]_{1 \le 
\beta \le s_j}]_{1\le j \le n_\infty}]$ is in $\ker \Gamma$ as claimed.

Conversely, if $[[c_{j \beta}]_{1 \le \beta \le s_j}]_{1 \le j \le n_\infty}$ 
is in $\ker \Gamma$ and we set
\[
h(p) = \sum_{j=1}^{n_\infty} \sum_{\beta =1}^{s_j} K(\widetilde{\chi}; p, 
\mu^j) c_{j \beta},
\]
then one can verify that $h \in {\cal O}(\widetilde{\chi} \otimes \Delta) 
({\boldsymbol \omega})(X)$ by reversing the steps of the above argument.
\end{pf}

The analogue of Lemma \ref{L1:absint} at the level of bundle 
endomorphisms is the following.

\begin{lemma}  \label{L1':absint}
Suppose that $T$ is a holomorphic bundle map from $\chi \otimes \Delta$ 
to $ (\widetilde{\chi} \otimes \Delta) ( {\boldsymbol \mu}, {\bold u})$ 
such that $T(q) = Q \colon 
\chi(q) \to \widetilde{\chi}(q)$.  Then there exists a unique choice of 
operators $\widehat{x}_{j \beta} \colon \widetilde{\chi}(q) \to {\bold C}$ 
such that
\[
T(p) = \left[ K(\widetilde{\chi}; p, q) +
\sum_{j=1}^{n_\infty} \sum_{\beta = 1} ^{s_j} K(\widetilde{\chi}; p, 
\mu^j) u_{j \beta} \widehat{x}_{j \beta} \right] Q K(\chi;p,q)^{-1}.
\]
\end{lemma}

\begin{pf}
Choose operators $\widehat{x}_{j \beta} \colon \widetilde{\chi}(q) \to 
{\bold C}$ so that
\[
\text{Res}_{\mu^j} T(\cdot) K(\chi; \mu^j,q)Q^{-1} = \sum_{ \beta}^{s_j} 
u_{j \beta} \widehat{x}_{j \beta}.
\]
Set
\[
\widehat{T}(p) = \left[ K(\widetilde{\chi}; p,q) + \sum_{j=1}^{n_\infty} 
\sum_{\beta = 1}^{s_j} K(\widetilde{\chi}; p, \mu^j) u_{j \beta} 
\widehat{x}_{j \beta} \right] Q(K(\chi, p,q)^{-1}.
\]
Then, for any vector $v \in \chi(q)$ we have
\begin{gather} \notag
T(\cdot) K(\chi; \cdot, q)v - \widehat{T}(\cdot) 
K(\chi; \cdot, q) v  \\ \notag
=T(\cdot) K(\chi; \cdot, q) v - \left[ K(\widetilde{\chi}; \cdot, q)Q
+\sum_{j=1}^{n_\infty} \sum_{\beta=1}^{s_j} K(\widetilde{\chi}; \cdot, 
\mu^j) u_{j \beta} \widehat{x}_{j \beta}Q \right] v \\
\end{gather}
is an element of ${\cal O}(\widetilde{\chi} \otimes \Delta)(X)$.  By our 
standing assumption that $h^0(\widetilde{\chi} \otimes \Delta) = 0$, we 
conclude that $(T-\widehat{T})(\cdot) K(\chi; \cdot, q)v=0$ for all $v 
\in \chi(q)$.  This is enough to force $T=\widehat{T}$, and the lemma follows.
\end{pf}

\begin{pf*}{Proof of Theorem \ref{T:absint}}

We first argue that  necessarily $\Gamma$ is square and invertible if 
a solution $T$ to the interpolation problem (ABSINT) exists.  
To do this, we show first  that $\ker \Gamma = 
\{0\}$ and secondly, that $\Gamma$ is square.

To see that $\ker \Gamma =\{0\}$, we proceed as follows.
If $T$ is a solution of (ABSINT), then multiplication by $T$ induces a 
biholomorphic bundle map between $\chi \otimes \Delta$ and 
$(\widetilde{\chi} \otimes \Delta)({\boldsymbol \omega})$.
By assumption, $h^0(\chi \otimes \Delta) = 0$.  Hence we also have 
$h^0( (\widetilde{\chi} \otimes \Delta) ({\boldsymbol  \omega})) = 0$.  Now it 
follows from Lemma \ref{L2:absint} that $\ker \Gamma = \{0\}$.

Next we argue that $\Gamma$ is square.  Since $\chi$ and 
$\widetilde{\chi}$ by assumption are both flat, both $\chi$ and 
$\widetilde{\chi}$ have  degree 0, as do $\det \chi$ and $\det 
\widetilde{\chi}$.  Then
\[
\text{deg}(\det T) = 
\text{deg}(\det \widetilde{\chi}) - \text{ deg} (\det \chi) = 0.
\]
 If $T$ is a solution of (ABSINT), then the total number of 
zeros $n_0(T)$ of $T$ (counted with multiplicities as appropriate 
for meromorphic matrix  functions---see Chapter 3 of \cite{bgr}) 
is equal to 
$\sum_{i=1}^{n_0} t_i =:N_0$ which is the number of rows of 
$\Gamma$, while the total number of poles $n_\infty(T)$ (again counted with 
multiplicities) is equal to $ \sum_{j=1}^{n_\infty} s_j =:N_\infty$ which 
is equal to
the number of columns of $\Gamma$.   In general we have $\text{deg}(\det T) 
= n_0(T) - n_\infty(T)$.  Hence the equality $\text{deg}(\det T) = 0$ for 
$T$ a solution
of (ABSINT) implies that $\Gamma$ is square.  Combining this with the 
result of the previous paragraph, we see that $\Gamma$ is invertible as well.

If $T$ is a solution of (ABSINT), then  in particular $T$ 
satisfies the hypotheses of Lemma \ref{L1':absint} and hence $T$ has the 
form
\begin{equation}  \label{ansatz}
T(p) = \left[ K(\widetilde{\chi};p,q) + \sum_{j=1}^{n_\infty} \sum_{\beta 
= 1}^{s_j} K(\widetilde{\chi}; p, \mu^j) u_{j \beta} \widehat{x}_{j 
\beta} \right] Q K(\chi, p,q)^{-1}
\end{equation}
for appropriate operators $\widehat{x}_{j \beta} \colon 
\widetilde{\chi}(q) \to {\bold C}$.  We now find what additional 
restrictions on $\widehat{x}_{j \beta}$ are forced by the zero and 
coupled zero-pole interpolation conditions (ii) and (iii) in (ABSINT).

Suppose that $i$ is an index for which $\lambda^i \ne \mu^j$ for 
any $j$.  Then the zero interpolation condition $x_{i \alpha}^T T(\lambda^i)
=0$ forces, for all $\alpha$ between $1$ and $t_i$,
\[
x_{i \alpha}^T K(\widetilde{\chi}; \lambda^i,q) Q K(\chi; \lambda^i, q)^{-1}
+ \sum_{j=1}^{n_\infty} \sum_{\beta=1}^{s_i}x_{i \alpha}^T 
 K(\widetilde{\chi}; \lambda^i, \mu^j) u_{j \beta}
  \widehat{x}_{j \beta} Q K(\chi; 
\lambda^i,q)^{-1} = 0.
\]
Recalling the definition of $\Gamma_{ij, \alpha \beta}$, we can rewrite 
this as
\begin{equation} \label{Gamma1}
\sum_{j=1}^{n_\infty} \sum_{\beta=1}^{s_j} \Gamma_{ij, \alpha \beta} \widehat{x}_{j \beta} = x_{i \alpha} 
K(\widetilde{\chi}, \lambda^i, q)
\end{equation}
for all index pairs $(i,\alpha)$ such that $\lambda^i \ne \mu^j$ for any $j$.

We next consider an index $i$ for which $\lambda^i = 
\mu^j:= \xi^{ij}$ for some $j$.  Let $x_{i \alpha}(p)$ be a local 
holomorphic section of $\widetilde{\chi}^\vee$ as in the third set of 
interpolation conditions.  Let $\varphi(p)$ be the meromorphic local 
section of $\chi$ given by
\[
\varphi(p) = \frac{K(\chi; p,q)}{\sqrt{dt^{ij}}(p)} e
\]
for a vector $e \in \chi(q)$ where $t^{ij}(p)$ is a local coordinate on 
$X$ centered at $\xi^{ij}$, and let $u(p)$ be the local meromorphic 
section of $\widetilde{\chi}$ given by $u(p) = T(p) \varphi(p)$.  From 
\eqref{ansatz} we have then
\[
u(p) = \frac{K(\widetilde{\chi};p,q)}{\sqrt{dt^{ij}}(p)} Q e
+ \sum_{j=1}^{n_\infty} \sum_{\beta = 1}^{s_j} 
\frac{K(\widetilde{\chi};p,\mu^j)}{\sqrt{dt^{ij}}(p)} u_{j \beta} 
\widehat{x}_{j \beta} Q e.
\]
Then the coefficients $[u]_{-1}$ and $[u]_0$ in the Laurent expansion of 
$u(p)$ centered at $\xi^{ij}$ with respect to local coordinate $t^{ij}$ 
are given by
\begin{align} \label{-1Laurent}
[u]_{-1} &= \sum_{\beta = 1}^{s_j} u_{j \beta} \sqrt{dt^{ij}}(\xi^{ij}) 
\widehat{x}_{j \beta} Q e, \\
\label{0Laurent}
[u]_0 &= \frac{K(\widetilde{\chi}; \xi^{ij},q)}{\sqrt{dt^{ij}}(\xi^{ij})}
Q e + \sum_{k \ne j} \sum_{\beta = 1}^{s_k} \frac{K(\widetilde{\chi}; 
\xi^{ij}, \mu^k)}{\sqrt{dt^{ij}}(\xi^{ij})} u_{k \beta}
 \widehat{x}_{k \beta} Q e \\
\notag
& + \sum_{\beta=1}^{t_j} \frac{A_\ell(\xi^{ij})}{dt(\xi^{ij})} 
\sqrt{dt}(\xi^{ij}) u_{j \beta} \widehat{x}_{j \beta} Q e
\end{align}
so the alternate coupled interpolation condition (iii) given by 
\eqref{coupledint} implies that  
\[
x_{i \alpha}^T \left( A(\xi^{ij}) [u]_{-1} + [u]_0 dt^{ij}(\xi^{ij}) \right)
= - \sum_{\beta = 1}^{s_j} \rho_{ij, \alpha \beta} \widehat{x}_{j \beta} Q e
\cdot \sqrt{dt^{ij}}(\xi^{ij}).
\]
Substitution of the expressions \eqref{-1Laurent} and \eqref{0Laurent} 
for $[u]_{-1}$ and $[u]_0$ gives
\begin{align} \notag
\sum_{\beta=1}^{s_j} x_{i \alpha}^T A(\xi^{ij}) u_{j \beta} 
\widehat{x}_{j \beta} Q e \sqrt{dt^{ij}}(\xi^{ij}) 
& + x_{i \alpha}^T 
\frac{K(\widetilde{\chi}; \xi^{ij}, q)}{\sqrt{dt^{ij}}(\xi^{ij})} Q e 
dt^{ij}(\xi^{ij}) \\
\notag
+ \sum_{k \ne j} \sum_{\beta = 1}^{s_k} x_{i \alpha}^T 
\frac{K(\widetilde{\chi}; \xi^{ij}, \mu^k)}{\sqrt{dt^{ij}}(\xi^{ij})}
u_{k \beta} \widehat{x}_{k \beta} Q e \cdot dt^{ij}(\xi^{ij}) 
& +
\sum_{\beta = 1}^{t_j} x_{i \alpha}^T 
\frac{A_\ell(\xi^{ij})}{\sqrt{dt^{ij}}(\xi^{ij})} u_{j \beta} 
\widehat{x}_{j \beta} Q e \cdot dt^{ij}(\xi^{ij}) \\
\notag
= &- \sum_{\beta = 1}^{s_j} \rho_{ij, \alpha \beta} \widehat{x}_{j \beta} 
Q e \cdot \sqrt{dt^{ij}}(\xi^{ij}).
\end{align}
By using the result of Lemma \ref{L:cauchyexp} and recalling the definition 
\eqref{defGamma} of $\Gamma_{ij, \alpha \beta}$ we see that this 
expression collapses to
\[
\sum_{k=1}^{n_\infty} \sum_{\beta = 1}^{s_k} \Gamma_{ik, \alpha \beta}
\widehat{x}_{k \beta} Q e  = x_{i \alpha}^T 
K(\widetilde{\chi}; \xi^{ij},q) Qe.
\]
Since this must hold for all $e \in \chi(q)$, we arrive at the operator 
equation
\begin{equation} \label{Gamma2}
\sum_{k=1}^{n_\infty} \sum_{\beta=1}^{s_k} \Gamma_{ik, \alpha \beta} 
\widehat{x}_{k \beta} = x_{i \alpha}^T K(\widetilde{\chi}; \xi^{ij},q)
\end{equation}
which must hold for all index pairs $(i, \alpha)$ for which $\lambda^i = 
\mu^j$ for some $j$. Combining \eqref{Gamma1} and \eqref{Gamma2} gives us
\[ 
\Gamma \widehat{{\bold x}} = K^{{\boldsymbol \lambda}, {\bold x}}(q)
\]
where we have set $\widehat{{\bold x}}$ equal to the column vector 
$[[ \widehat{x}_{j \beta}]_{ 1 \le \beta \le s_j}]_{1 \le j \le n_\infty}$. 
 Plugging this 
value into \eqref{ansatz} leaves us with the formula \eqref{solution} for 
the solution $T$.  This also establishes the uniqueness of the solution 
of (ABSINT) whenever it exists.

Since $T$ is analytic at the points $p^1, \dots, p^n$ where $K(\chi; 
\cdot, q)^{-1}$ has poles, necessarily the residue conditions 
\eqref{residues} must hold as well.  The necessity and uniqueness parts 
of the theorem are now established.

Conversely, assume that $\Gamma$ is invertible and that the residue 
conditions \eqref{residues} hold.  We define $T(p)$ by the formula 
\eqref{solution}.  Then $T$ is a meromorphic bundle map of $\chi$ and 
$\widetilde{\chi}$ with only simple poles which occur at most at the 
points $\mu^1, \dots, \mu^{n_\infty}$ with 
\[
\text{im Res}_{\mu^i} T(\cdot) \subset \text{span }\{u_{j \alpha} \colon 
1 \le \alpha \le s_j\}
\]
for $j=1, \dots, n_\infty$.  Since $T$ has the form \eqref{ansatz} with 
operators $\widehat{x}_{j \beta}$ ($1 \le j \le n_\infty$, $1 \le \beta 
\le s_j$) satisfying \eqref{Gamma1} and \eqref{Gamma2}, we see that also 
$T$ satisfies the interpolation conditions (ii) and (iii) in (ABSINT) 
as well.  Hence, the number of poles $n_\infty(T)$ of $T$ (counting 
multiplicities for a meromorphic matrix function as in Chapter 3 of 
\cite{bgr}) is at most 
$\sum_{j=1}^{n_\infty} s_j =: N_\infty$ and the number of zeros $n_0(T)$ 
of $T$ (counting multiplicities) is at least 
$\sum_{i=1}^{n_\infty} t_i =: N_0$.  As $T$ is a bundle map of the flat 
bundles $\chi$ and $\widetilde{\chi}$, we know that $n_0(T) = 
n_\infty(T)$.  On the other hand, since $\Gamma$ is square we have 
$N_0=N_\infty$.  From the chain of inequalities
\[
N_0 \le n_0(T) = n_\infty(T) \le N_\infty
\]
combined with the equality $N_0 = N_\infty$, we get that $n_0(T) = N_0$ 
and $n_\infty(T) = N_\infty$.  This implies that necessarily
\[
\text{im Res}_{p = \mu^j} T(p) = \text{span}\{u_{j1}, \dots, u_{j s_j} \}
\]
and
\[
{\text{im Res}}_{p=\lambda^i}(T^{\vee})^{-1}(p) = \text{span}\{x_{i1}, \dots, 
x_{i t_i}\}
\]
and that $T(p)$ is analytic and invertible at every point $p \in X$ 
outside of $\mu^1, \dots, \mu^{n_\infty}$, $\lambda^1, \dots, \lambda^{n_0}$.
This verifies that $T$ is a bona fide solution of the interpolation 
problem (ABSINT).

It remains only to verify the formula \eqref{inversesolution} for the 
inverse of $T$.  To see this, we note that $(T^{-1})^T$ is also the 
solution of an interpolation problem of the type (ABSINT), namely the one 
with data set $[{\boldsymbol \omega}]^\vee$ given by 
\begin{enumerate}
\item $( {\boldsymbol \lambda},{\bold x})$ in place of 
$( {\boldsymbol \mu}, {\bold u})$,
\item $({\boldsymbol \mu}, {\bold u})$ in place of $({\boldsymbol \lambda},
 {\bold x})$, and
\item $-\rho_{ji, \beta \alpha}$ in place of $\rho_{ij, \alpha \beta}$.
\end{enumerate}
The matrix $\Gamma$ associated with this interpolation problem turns out 
to be exactly $-\Gamma^T$ where $\Gamma$ is as in \eqref{defGamma}.  
Hence by Theorem \ref{T:absint} the bundle map $(T^{-1})^T$ must be given by
\begin{equation} \label{Tinvtr}
(T^{-1})^T = [K(\widetilde{\chi}^\vee;p,q) - K_{{\boldsymbol \lambda}, 
{\bold x}}(p) (\Gamma^{-1})^T K^{{\bold u}, 
{\boldsymbol \mu}}(q)] (T(q)^{-1})^T K(\chi^\vee; 
p,q)^{-1}.
\end{equation}
By the uniqueness property of Cauchy kernels, it is easy to see that
\[
K(\chi^\vee; p,q)^T = - K(\chi; q,p).
\]
Hence,  taking transpose on both sides of \eqref{Tinvtr} gives
\begin{align} \notag
T^{-1}(p) = & -K(\chi;q,p)^{-1} T(q)^{-1}[-K(\widetilde{\chi}; q,p) 
- K_{ {\boldsymbol \mu}, {\bold u}}(q) \Gamma^{-1} K^{{\bold x},
{\boldsymbol \lambda} }(p) ] \\
\notag = & K(\chi;q,p)^{-1} T(q)^{-1}[-K(\widetilde{\chi}; q,p) 
+ K_{{\boldsymbol \mu}, {\bold u}}(q) \Gamma^{-1} K^{{\bold x},{\boldsymbol
 \lambda}}(p) ]
\end{align}
and the formula \eqref{inversesolution} follows.
\end{pf*}

{\bf Remark:} In case $\widetilde{\chi} = \chi =:\chi_0$ are both taken 
to be the trivial bundle of rank $r$, the Cauchy kernel $K(\chi_0; \cdot, 
\cdot)$ has the scalar form $k_0(\cdot, \cdot) I_r$ where $k_0(\cdot, 
\cdot)$ is the Cauchy kernel for the trivial line bundle over $X$.  In 
this case the ansatz \eqref{ansatz} simplifies to
\[
T(p) = Q + \sum_{j=1}^{n_\infty} \sum_{\beta=1}^{s_j} f_{\mu^j}(p) u_{j 
\beta} \widehat{x}^\prime_{j \beta}
\]
where we have set
\[f_{\mu}(p) = \frac{k_0(p,\mu)}{k_0(p,q)}
\]
and the row vectors $\widehat{x}^\prime_{j \beta} = \widehat{x}_{j \beta} 
Q$ are now taken to be the unknowns.  Note that $k_0(\cdot, q)$ is a 
half-order differential with divisor of degree $g-1$ and a pole at $q$; 
if we assume that the zeros are distinct, this divisor has the form $p^1 + 
\dots + p^g - q$ for distinct points $p^1, \dots, p^g,q \in X$.  If the 
image of the divisor $p^1 + \dots + p^g$ under the Abel-Jacobi map is not 
on the classical theta divisor in the Jacobian (i.e. if $(p^1 + \dots + 
p^g$ is a {\it non-special} divisor), then there are no nonzero constant
meromorphic functions with only poles equal to at most simple poles at
 the points $p^1, \dots, p^g$; this 
corresponds to our assumption that $h^0(\chi_0 \otimes \Delta)=h^0(\Delta)=0$.  
Furthermore, in this case, the global scalar meromorphic function 
$f_\mu(\cdot)$ on $X$ (for $\mu$ a point of $X$ disjoint from $p^1, 
\dots, p^g,q$) is uniquely determined (up to a nonzero scalar multiple) 
by the condition that it have a pole at $\mu$ and that its divisor 
$(f_\mu)$ satisfy 
\[
(f_\mu) \ge q - \mu - p^1 - \dots - p^g. 
\]
In this way our results and analysis on the (ABSINT) problem reduce to 
the work in \cite{bc} for the trivial bundle case.
Notice that $\chi_0$ can be replaced here by $\xi\otimes\chi_0$
for any line bundle $\xi$ of degree $0$ satisfying $h^0(\xi \otimes \Delta)=0$,
replacing the Cauchy kernel $k_0(\cdot,\cdot)$ for the trivial line bundle
by the Cauchy kernel for $\xi$;
this corresponds to letting $p^1 + \dots + p^g$ be any non-special
effective divisor of degree $g$.
\vspace{.1in}

One remaining piece of business in this section is the proof of Theorem 
\ref{T:bvabsint}.  The problem of identifying the unknown input bundle in 
a more explicit fashion will be addressed in Section \ref{S:conint}.

\begin{pf*}{Proof of Theorem \ref{T:bvabsint}}

If there exists such an input bundle $\chi$ and meromorphic bundle map 
$T$, then $T$ implements a biholomorphic bundle map between $\chi \otimes 
\Delta$ and $(\widetilde{\chi} \otimes \Delta)({\boldsymbol \omega})$.  
Since $h^0(\chi 
\otimes \Delta) = 0$, it then follows 
from Lemma \ref{L2:absint} that $\Gamma$ is injective.  Since deg$(\chi 
\otimes \Delta) = r(g-1)$, it must be the case that deg$(\widetilde{\chi} 
\otimes \Delta)({\boldsymbol \omega})=r(g-1)$ as well.  This means that 
$\Gamma$ is square, 
and hence invertible.

Conversely, suppose that $\Gamma$ is square and invertible.  Define a 
bundle $\chi$ so that $\chi \otimes \Delta \cong (\widetilde{\chi} 
\otimes \Delta)({\boldsymbol \omega})$.  Since $\Gamma$ is square, 
it follows that
\[ 
\text{deg}\left( (\widetilde{\chi}\otimes \Delta)({\boldsymbol \omega})\right) = 
\text{deg}(\widetilde{\chi} \otimes \Delta) = r(g-1),
\]
and hence $\text{deg}(\chi \otimes \Delta) = r(g-1)$.  Since $\ker \Gamma 
= \{0\}$, we know by Lemma \ref{L2:absint} that $h^0( (\widetilde{\chi} 
\otimes \Delta)({\boldsymbol \omega}))=0$; thus $h^0(\chi \otimes \Delta) = 0$.  It 
follows from these two facts as in the proof of Theorem 3.1 in \cite{hip} 
that $\chi$ is flat.

Let now $S \colon \chi \otimes \Delta \to (\widetilde{\chi} \otimes 
\Delta)({\boldsymbol \omega})$ be an implementation of the holomorphic bundle 
isomorphism between $\chi \otimes \Delta$ and $(\widetilde{\chi} \otimes 
\Delta)({\boldsymbol \omega})$.  
Define $T \colon \chi \to \widetilde{\chi}$ so that $S 
= T \otimes I_{{\cal O}(\Delta)}$.  Then $T$ is a meromorphic bundle map 
from $\chi$ to $\widetilde{\chi}$ which solves the interpolation problem 
(ABSINT).

\end{pf*}

\section{The line bundle case and Fay's identity}
\label{S:linebundle}

In this section we specialize the work of the preceding sections to the 
line bundle case.

We shall need here some basic facts concerning the Jacobian variety, 
the Abel-Jacobi map and associated theta functions (theta function, theta 
functions with characteristics and prime form)
for the Riemann surface $X$.  The review here is quite 
sketchy; for complete details the reader should consult \cite{oldfay}, 
\cite{farkaskra} or \cite{mumford}.
 
When the rank $r$ of the vector bundle $\chi$ is 1, one can get an 
explicit formula for $K(\chi; \cdot, \cdot)$ in terms of the Abel-Jacobi 
map for the surface $X$ and various variants of the classical theta 
function associated with the Jacobian variety of $X$ (see \cite{hip}).
Specifically, in this case we may assume that $\chi$ is a flat unitary line 
bundle with factor of automorphy (also called $\chi$) given by $\chi(A_j) 
= \exp (- 2 \pi i a_j)$, $\chi(B_j) = \exp (2 \pi i b_j)$, $j=1, \dots, g$, 
where $A_1, \dots, A_g, B_1, \dots, B_g$ form a canonical integral 
homology basis on $X$.  Let $\Omega$ be the corresponding period matrix, 
let $J(X)= {\bold C}^g / {\bold Z}^g + \Omega {\bold Z}^g$ be the 
Jacobian variety of $X$ and let $\phi \colon X \to J(X)$ be the 
Abel-Jacobi map.  As is standard, we extend $\phi$ by linearity to any 
divisor on $X$, and, using the correspondence between linear equivalence 
classes of divisors and isomorphism classes of line bundles, we consider 
$\phi$ to be defined on any line bundle on $X$ as well.  One can verify 
that then $\phi(\chi) = z$ where $z = \Omega a + b$ and $a,b \in {\bold 
R}^g$ have respective coordinates $a_j,b_j$. 
 Then the explicit formula 
for the Cauchy kernel (as given in \cite{hip}) is the following.  The 
verification is straightforward, once one has in hand the properties and 
factors of automorphy for the various objects involved.

\begin{theorem} \label{T:scalarCauchykernel}
For the case where $\chi$ is a flat unitary line bundle as above, the Cauchy 
kernel as defined in Section \ref{S:cauchyker} is given explicitly by 
\begin{equation}
K(\chi;p,q) = \frac{ \theta \begin{bmatrix} a \\ b \end{bmatrix} (\phi(q) 
- \phi(p)) } {\theta \begin{bmatrix} a \\ b \end{bmatrix} (0) E(q,p)}.
\label{lineCauchyker}
\end{equation}
\end{theorem}

In the statement of Theorem \ref{T:scalarCauchykernel}
$\theta \begin{bmatrix} a \\ b \end{bmatrix} ( \cdot)$ is the 
associated theta function with characteristics $\begin{bmatrix} a \\ b 
\end{bmatrix}$, $E(\cdot, \cdot)$ is the prime form on $X \times X$, and 
we assume the line bundle $\Delta$ of differentials of order 
$\frac{1}{2}$ has been chosen so that $\phi(\Delta) = - \kappa$, where 
$\kappa \in J(X)$ is Riemann's constant (see \cite{oldfay} and 
\cite{mumford}).  Note that a consequence of Riemann's theorem is that 
$\theta(z) \ne 0$ if and only if $h^0(\chi \otimes \Delta) = 0$, and 
hence $\theta \begin{bmatrix} a \\ b \end{bmatrix} (0) \ne 0$ in 
\eqref{lineCauchyker} and the formula makes sense.  No such explicit 
formula is known at present for the higher rank case except in genus 1
(see \cite{bcv}).

In the line bundle case one can also give an explicit formula for the 
canonical connections $\nabla_\chi$, $\nabla_\chi^*$ associated with the 
flat unitary line bundle $\chi$.  This is the content of the following 
Proposition.

\begin{proposition} \label{P:connection}
For the case where $\chi$ is a flat unitary line bundle with 
normalizations as above, then the canonical connections $\nabla_\chi$ 
and $\nabla_\chi^*$ are given by
 \begin{align}
\nabla_\chi y = & \left[ \sum_{j=1}^g \frac{\partial}{\partial z_j} \log 
\theta \begin{bmatrix} a \\ b \end{bmatrix} (0) \ \omega_j(p) \right] y + 
dy \\
= & \left[ \sum_{j=1}^g [2 \pi i a_j + \frac{\partial}{\partial z_j} \log 
\theta(z) ] \omega_j(p) \right] y + dy,\\
\nabla_\chi x = & - \left[ \sum_{j=1}^g \frac{\partial}{\partial z_j} \log 
\theta \begin{bmatrix} a \\ b \end{bmatrix} (0) \ \omega_j(p) \right] x + 
dx \\
= & - \left[ \sum_{j=1}^g [2 \pi i a_j + \frac{\partial}{\partial z_j} \log 
\theta(z) ] \omega_j(p) \right] x + dx.
\end{align}
\end{proposition}

\begin{pf} 
This follows directly by comparing the general expansion for the Cauchy kernel
$$
\frac{K(\chi; p,p_0)}{\sqrt{dt}(p) \sqrt{dt}(p_0)} =
\frac{1}{t(p) - t(p_0)} \left[ I_r + \frac{A_\ell}{dt}(p_0) t(p) + O(|t(p)|^2)
\right]
$$
on the one hand and substituting the expansion of the theta function
$$
\theta \begin{bmatrix} a \\ b \end{bmatrix} \left(\phi(p_0) - \phi(p)\right)
= \theta \begin{bmatrix} a \\ b \end{bmatrix} (0) - \sum_{j=1}^g 
\frac{\partial \theta \begin{bmatrix} a \\ b \end{bmatrix}} {\partial z_j} 
( 0 ) \frac{\omega_j}{dt} \left( p_0 \right) t(p) + O(|t(p)|^2)
$$
and the expansion of the prime form (see Corollary 2.5 in \cite{oldfay})
$$
E(p_0,p) = t(p) + O(|t(p)|^3)
$$
into \eqref{lineCauchyker}.
\end{pf}

{\bf Remark.}  From the formula for $\nabla_\chi$ and $\nabla_\chi^*$ it 
follows that the coefficients $A(p)$ and $A_\ell(p)$ are independent of 
the choice of homology bases (i.e., marking) on the Riemann surface $X$
as long as the bundle $\Delta$ of half-order differentials defined by 
$\phi(\Delta) = \kappa$ remains the same, 
since the unitary flat representative for a flat line bundle is unique.
It is an amusing exercise to verify this independence directly by using 
the transformation law for theta functions (see \cite{mumford} and 
\cite{Igusa}).

We next specialize the work of Section \ref{S:absint} to the scalar (or 
line bundle) case, where $\chi$ and $\widetilde{\chi}$ are flat unitary 
line bundles.  As explained in Section \ref{S:absint}, necessarily the 
multiplicities $s_j$ and $t_i$ are all 1 and without loss of generality 
we may take $u_{j1}=1$, $x_{i1}=1$ for all $i$ and $j$.  
Then the 
compatibility condition \eqref{comp} forces the third interpolation 
condition to be absent.  The data of the problem consists simply of the 
set of $n_\infty + n_0$ distinct points $\mu^1, \dots, \mu^{n_\infty}, 
\lambda^1, \dots, \lambda^{n_0}$ together with the flat unitary line bundles
$\chi$ and $\widetilde{\chi}$.  
The problem then is to produce a bundle 
map $T \colon \chi \to \widetilde{\chi}$ with divisor equal to 
${\boldsymbol{\lambda}} - 
{\boldsymbol{\mu}}$ (where we have set ${\boldsymbol{\lambda}} = \lambda^1 
+ \dots + \lambda^{n_0}$ and ${\boldsymbol \mu} =  
\mu^1 + \dots + \mu^{n_\infty}$).
If we view the bundles in terms of factors of automorphy, we can view $T$ 
simply as a multivalued function on $X$ having divisor equal to 
${\boldsymbol \lambda} - {\boldsymbol \mu}$ and factor of automorphy 
$\chi_T$ given by 
\[
\chi_T(A_j) = e^{-2 \pi i a_j}, \ \chi_T(B_j) = e^{2 \pi i b_j}
\text{ for } j=1, \dots, g
\]
where $\phi(\widetilde{\chi}) - \phi(\chi) = \Omega a + b$ (where we have 
set $a= \begin{bmatrix} a_1 & \dots & a_g \end{bmatrix}^T$ and $b= 
\begin{bmatrix} b_1 & \dots & b_g \end{bmatrix}^T$).

In the genus zero case where $X= \bold{C} \cup \{\infty\}$ is the Riemann 
sphere, any flat unitary line bundle is trivial and the problem is to 
produce a global meromorphic function with divisor equal to ${\boldsymbol
\lambda}
 - {\boldsymbol \mu}$. 
Trivially a solution exists if and only if $n_0=n_\infty$ and then the 
unique solution with value 1 at infinity is given in the multiplicative 
form
\begin{equation}
T(z) = \dfrac{\prod_{i=1}^{n_0} (z-\lambda^i)}{\prod_{j=1}^{n_\infty} 
(z-\mu^j)}. 
\label{g0prod}
\end{equation}
or in the partial fraction form
\begin{equation}
T(z) = 1 + \sum_{j=1}^{n_\infty} c_j (z-\mu^j)^{-1}
\label{g0partialfrac}
\end{equation}
where $c^T= \begin{bmatrix} c_1 & \dots & c_{n_\infty} \end{bmatrix}$ is 
the unique solution of the linear system of equations $S c 
= \begin{bmatrix} 1 & 
\dots & 1 \end{bmatrix}^T$ with $S$ equal to the Sylvester matrix
\[
 S = [S_{ij}] \text{ with } S_{ij} = \frac{1}{\mu^j - \lambda^i},
 \]
or, in other words, 
\begin{equation}
c = S^{-1} \begin{bmatrix} 1 & \dots & 1 \end{bmatrix}^T.
\label{g0partialfraccoef}
\end{equation}
It is possible to evaluate the vector $c$  explicitly from 
\eqref{g0partialfraccoef} once one knows the entries of $S^{-1}$ 
explicitly.  This in turn can be done once one knows an explicit formula 
for the determinant of a Sylvester matrix $S$, since the cofactor 
matrices are again of the same form.  In this way one can verify directly 
the equivalence of the two formulas \eqref{g0prod} and 
\eqref{g0partialfrac}.  For details on the algebra of 
this computation, we refer to 
Theorem 4.3.2 of \cite{bgr}.

We shall see that an analogous pair of formulas holds for the solution of
the abstract interpolation problem 
(ABSINT) for the higher genus case (for the line bundle setting).  This is 
the content of the next Theorem.  

\begin{theorem} \label{T:lineabsint}
Consider the problem (ABSINT) for the case where $\chi$ and 
$\widetilde{\chi}$ are flat unitary line bundles and given data set equal to 
$\boldsymbol{\lambda} - \boldsymbol{\mu} = \lambda^1 + \dots + \lambda^{n_0} - 
\mu^1 - \dots - 
\mu^{n_\infty}$ as above.  Then a solution exists if and only if
\begin{equation}
n_0 = n_\infty \text{ and } 
\phi(\widetilde{\chi}) - \phi(\chi) = \phi({\boldsymbol \lambda}) -
 \phi(\boldsymbol \mu).
\label{necessity}
\end{equation}
In this case,  if $q$ is a point of $X$ disjoint from the set 
$\{\lambda^1, \dots, \lambda^{n_0}, \mu^1, \dots, \mu^{n_\infty}\}$ of 
prescribed zeros and poles and $Q$ is any invertible fiber map from 
$\chi(q)$ onto $\widetilde{\chi}(q)$, then a solution of (ABSINT) having 
value $Q$ at $q$ is given in multiplicative form as
\begin{equation} T(p) =  
\dfrac{ \prod_{i=1}^{n_0}   
E(p,\lambda^i)/E(q,\lambda^i)}{\prod_{j=1}^{n_\infty} E(p,\mu^j)/E(q,\mu^j)}
\exp(-2 \pi i a^T (\phi(p) - \phi(q)) Q
\label{prod}
\end{equation}
where $a^T = \begin{bmatrix} a_1 & \dots & a_g \end{bmatrix}$ and 
$\phi({\boldsymbol{\lambda}}) - \phi({\boldsymbol{\mu}}) = 
\Omega a + b$ with $a,b \in {\bold R}^g$.

If $\chi \otimes \Delta$ and $\widetilde{\chi} \otimes \Delta$ have no 
nontrivial holomorphic sections, then the solution 
$T$ with $T(q)=Q$ is unique and
alternatively is given by the partial fraction formula
\begin{align}
T(p)=& \left\{ \dfrac{\theta[\widetilde{z}](\phi(q) - 
\phi(p))}{\theta[\widetilde{z}](0) E(q,p)}
+ \sum_{j=1}^{n_\infty} \sum_{k=1}^{n_0} 
\dfrac{\theta[\widetilde{z}](\phi(\mu^j) - 
\phi(p))}{\theta[\widetilde{z}](0)E(\mu^j,p)} \cdot  [\Gamma^{-1}]_{jk}
\cdot \dfrac{\theta[\widetilde{z}](\phi(q)- 
\phi(\lambda^k))}{\theta[\widetilde{z}](0) E(q,\lambda^k) } \right\} \notag \\
& \times Q \cdot \dfrac{\theta[z](0) E(q,p)}{\theta[z](\phi(q) - \phi(p))}
\label{partialfrac}
\end{align}
where the $n_0 \times n_\infty$ matrix $\Gamma$ is given by
\[
\Gamma = [\Gamma_{ij}] \text{ with }
\Gamma_{ij} =- \dfrac{\theta[\widetilde{z}](\phi(\mu^j) - 
\phi(\lambda^i))}{\theta[\widetilde{z}](0) E(\mu^j, \lambda^i)}.
\]
and where we have set
\[
z = \phi(\chi), \ \widetilde{z} = \phi(\widetilde{\chi}).
\]
Here we write $\theta[z](\lambda)$ rather than $\theta\begin{bmatrix} a 
\\ b \end{bmatrix}(\lambda)$ if $z = \Omega a + b \in \bold{C}^g$ with 
$a, b \in {\bold R}^g$.  
\end{theorem}

\begin{pf} If such a bundle map exists, then the bundle $\chi \otimes 
{\cal O}(\boldsymbol{\lambda} - \boldsymbol{\mu})$ and 
$\widetilde{\chi}$ are holomorphically equivalent.  In particular, the 
divisor $\boldsymbol{\lambda} - \boldsymbol{\mu}$ 
must have degree 0 since both $\chi$ and 
$\widetilde{\chi}$ are flat bundles.  The equality $\phi(\widetilde{\chi}) =
\phi(\chi) + \phi(\boldsymbol{\lambda}) - \phi(\boldsymbol{\mu})$ 
then follows from the 
correspondence between flat bundles and linear equivalence classes of 
divisors mentioned above. This verifies the necessity condition 
\eqref{necessity}. 

Conversely,  assume that \eqref{necessity} holds and define $T$ by the 
right-hand side of \eqref{prod}. 
That the zero-pole divisor of $T$ 
is $\boldsymbol{\lambda} - \boldsymbol{\mu}$ follows directly from the fact 
that the divisor of the 
prime form $(p,q) \to E(p,q)$ is the diagonal $\{(p,p) \colon p \in X\} 
\subset X \times X$.  One can next check from the known period relations 
of $E$ that the right-hand side of \eqref{prod} has the factor of 
automorphy $\chi_T$
\[
\chi_T(A_j) = \exp (-2 \pi i a_j), \ \chi_T(B_j) = \exp (2 \pi i b_j) 
\text{ for } j=1, \dots, g
\]
where $a_1, \dots, a_g, b_1, \dots, b_g$ are respective components of
$a,b \in {\bold R}^g$ chosen so that $\Omega a + b = 
\phi(\boldsymbol{\lambda}) - 
\phi(\boldsymbol{\mu})$.  The second condition in \eqref{necessity} now 
guarantees that $T(\cdot)$ so defined is a 
bundle map from $\chi$ into $\widetilde{\chi}$. The uniqueness assertion 
is a consequence of Lemma \ref{L1':absint}.

The alternative formula \eqref{partialfrac} is  
 simply a rewriting of the formula \eqref{solution} from Theorem 
\ref{T:absint} specialized to the line bundle case, 
where we have substituted the explicit formula 
\eqref{lineCauchyker} for the Cauchy kernel from 
Theorem \ref{T:scalarCauchykernel}.

\end{pf}

Note that part of the content of Theorem \ref{T:lineabsint} is that the 
matrix $\Gamma$ is invertible whenever $\chi$ and $\widetilde{\chi}$ have 
no nontrivial holomorphic sections and (ABSINT) has a solution.

 It is of interest to specialize Theorem \ref{T:lineabsint} to the case 
 of one prescribed zero and pole $\lambda - \mu = \lambda^1 - \mu^1$.
 When this is done we obtain the following result.
 
 \begin{theorem} \label{T:redlineabsint}  
 If $\chi$ and $\widetilde{\chi}$ are two flat unitary line bundles 
 such that neither $\chi \otimes \Delta$ nor $\widetilde{\chi}
  \otimes \Delta$ have 
 nontrivial holomorphic sections and $\lambda$ and $\mu$ are two distinct 
 points of $X$, then the unique meromorphic bundle map from $\chi$ to 
 $\widetilde{\chi}$ with zero-pole divisor equal to $\lambda - \mu$  
 and value $Q \neq 0$ at the point $q \in X$ is given by either
 \begin{equation} \label{specialprod}
 T(p) = \dfrac{E(p,\lambda)}{E(p, \mu)} \dfrac{E(q,\mu)}{E(q,\lambda)}
      \exp(-2 \pi i a^T(\phi(p) - \phi(q)) Q
      \end{equation}
 or 
 \begin{align}
 T(p) =&  \exp(-2 \pi i a^T(\phi(p) - \phi(q)))
 \left\{ \dfrac{\theta(z + \phi(\lambda) - \phi(\mu) + \phi(q) - 
 \phi(p)) \theta(z) }{\theta(z + \phi(\lambda) - \phi(\mu)) 
 \theta(z + \phi(q) - \phi(p))} \right.  \notag \\
 & \left. - \dfrac{\theta(z + \phi(\lambda) - \phi(p)) \theta(z+ \phi(q) - 
 \phi(\mu)) E(\mu, \lambda) E(q,p)}
 {\theta(z+\phi(\lambda) - \phi(\mu)) \theta(z + \phi(q) - \phi(p)) 
 E(\mu,p) E(q,\lambda)} \right\} Q.
 \label{specialpartialfrac}
 \end{align}
 \end{theorem}
 
 \begin{pf}  The starting point of course is formula \eqref{prod} and 
 \eqref{partialfrac} specialized to the case $\boldsymbol{\lambda} - 
 \boldsymbol{\mu} 
 = \lambda - \mu$. The formula \eqref{specialprod} is an immediate 
 consequence of \eqref{prod}.  Derivation of \eqref{specialpartialfrac}
 requires a little bit of algebra.  We use the definition 
 \[
 \theta[z](\lambda) = \exp( \pi i a^T \Omega a + 2 \pi i a (\lambda + b)) 
 \theta(\lambda + z)
 \]
 (where $z = \Omega a + b$ with $a, b \in {\bold R}^g$) to express all 
 theta functions with characteristic $\theta[z](\cdot)$ in terms of the 
 theta function itself $\theta(\cdot)$.  When this is plugged into 
 \eqref{partialfrac} and little bit of algebra is used to collect the 
 exponential factor (noting that $a = a_{\widetilde{z}} - a_z$
 if $\widetilde{z} - z (=\phi(\lambda) - \phi(\mu)) = \Omega a + b$ and
 $\widetilde{z} = \Omega a_{\widetilde{z}} + b_{\widetilde{z}}$, 
 $z = \Omega a_z + b_z$ with $a$, $b$, $a_{\widetilde{z}}$, 
 $b_{\widetilde{z}}$, $a_z$, $b_z$ in ${\bold R}^g$), we get
 \begin{align}
 T(p) = & \exp(- 2 \pi i a^T(\phi(p) - \phi(q) ) ) \left\{
 \dfrac{ \theta(z + \phi(\lambda) - \phi(\mu) + \phi(q) - 
 \phi(p))}{\theta(z + \phi(\lambda) - \phi(\mu)) E(q,p)} \right. \notag \\
 &  - \dfrac{\theta(z + \phi(\lambda) - 
 \phi(p))}{\theta(z+\phi(\lambda) - \phi(\mu)) E(\mu,p)} \cdot
 \dfrac{ \theta(z + \phi(\lambda) - \phi(\mu) E(\mu, \lambda)}{\theta(z)}
 \notag \\
  & \left. \cdot \dfrac{\theta(z+\phi(q) - \phi(\mu))}
  {\theta(z + \phi(\lambda) - \phi(\mu)) E(q, \lambda)} \right\} 
  \cdot Q \dfrac{\theta(z) E(q,p)}{\theta(z+ \phi(q) - \phi(p))}. 
 \notag
 \end{align}
 The formula \eqref{specialpartialfrac} now follows by simple algebraic 
 manipulation.
 \end{pf}
 
 As a Corollary we obtain a version of Fay's Trisecant Identity (see 
 \cite{oldfay} formula (45) page 34 or \cite{mumford} Volume II page 3.214).
 
 \begin{corollary} \label{C:faytrisecant}
 For $X$ a compact Riemann surface, $\phi$ its Abel-Jacobi map, 
 $\theta(\lambda)$ and $E(p,q)$ its associated respective 
 theta function and prime 
 form, $p,q,\lambda, \mu$ points of $X$ and $z \in {\bold C}^g$, the 
 following identity holds:
 \begin{gather}
 \theta(z+\phi(\lambda) - \phi(\mu)) \theta(z+\phi(q) - \phi(p)) 
 E(p,\lambda) E(q,\mu)  \notag \\
  + \theta(z+\phi(\lambda) - \phi(p)) \theta(z + 
 \phi(q) - \phi(\mu)) E(\lambda, \mu) E(q,p) \notag \\
 =  \theta(z + \phi(\lambda) - \phi(\mu) + \phi(q) - \phi(p)) \theta(z) 
 E(p,\mu) E(q,\lambda).
 \label{faytrisecant}
 \end{gather}
 
 \end{corollary}
 
 \begin{pf} From the identity of the two expressions \eqref{specialprod} 
 and \eqref{specialpartialfrac} for $T(p)$ in Theorem 
 \ref{T:redlineabsint}, we have equality of the following two expressions 
 for $\exp(2 \pi i a^T(\phi(p) - \phi(q))) T(p) Q^{-1}$:
 \begin{align}
 \dfrac{E(p,\lambda)}{E(p,\mu)} \dfrac{E(q,\mu)}{E(q,\lambda)} = &
   \dfrac{\theta(z + \phi(\lambda) - \phi(\mu) + \phi(q) - 
 \phi(p)) \theta(z) }{\theta(z + \phi(\lambda) - \phi(\mu)) 
 \theta(z + \phi(q) - \phi(p))}   \notag \\
 &  - \dfrac{\theta(z + \phi(\lambda) - \phi(p)) \theta(z+ \phi(q) - 
 \phi(\mu)) E(\mu, \lambda) E(q,p)}
 {\theta(z+\phi(\lambda) - \phi(\mu)) \theta(z + \phi(q) - \phi(p)) 
 E(\mu,p) E(q,\lambda)} 
 \notag
 \end{align}
 Multiplication of both sides by $\theta(z+\phi(\lambda)-\phi(\mu)) 
 \theta(z + \phi(q) - \phi(p)) E(\mu,p) E(q, \lambda)$ along with
  a liberal use of 
 the general identity $E(p,q) = -E(q,p)$ along with some algebra now 
 leads to Fay's identity \eqref{faytrisecant} as desired.
 \end{pf}
 
 The identity \eqref{faytrisecant} is actually a special case of a more 
 general identity (see Corollary 2.19 in \cite{oldfay}) which gives an 
 explicit expression for the determinant of a matrix $M$ of the form
 \[
 M = [M_{ij}] \text{ where } M_{ij} = \dfrac{\theta(z + \phi(\mu^j) 
 - \phi(\lambda^i))}{E(\mu^j, \lambda^i)}.
 \]
 Since the cofactor matrices of such a matrix are of the same form, one 
 can then compute explicitly (in terms of theta functions and prime 
 forms) the entries of the inverse of the matrix $\Gamma$ appearing in 
 Theorem \ref{T:lineabsint}.  In this way one can verify by direct 
 computation the identity of the two expressions \eqref{prod} and 
 \eqref{partialfrac} for $T(p)$ in Theorem \ref{T:lineabsint}.  This then 
 is a canonical higher genus generalization of Theorem 4.3.2 in \cite{bgr}.

Note that this proof of Fay's identity arises from equating a 
multiplicative formula for the solution of a zero-pole interpolation 
problem to a partial-fraction expression for the same solution.  Formula 
\eqref{solution} in Theorem 
\ref{T:absint} gives an analogue of the partial fraction expression for 
the solution of a zero-pole interpolation problem for a vector bundle 
endomorphism.  Formula \eqref{solution}, giving a connection between the Cauchy 
kernel $K(\chi; p,q)$ and $K(\widetilde \chi; p,q)$, can be viewed as a 
matrix-valued version of the Fay trisecant identity.

There is one case in higher rank when a multiplicative representation does 
exist, namely the case of full rank zero-pole interpolation (where the 
given pole vectors $\{u_{j \beta} \colon 1 \le \beta \le s_j = r\}$ span 
the fiber space $\widetilde{\chi}(\mu^j)$ and 
the given null vectors $\{ x_{i \alpha} \colon 1 \le \alpha \le t_i =r\}$ 
span the fiber space $\widetilde{\chi}^\vee(\lambda^i)$) for each $i$ and 
$j$.  Then $T(p)$ is again given by \eqref{prod} where $Q$ now is a 
product of a scalar from the fiber of ${\cal O}({\boldsymbol\lambda} - 
{\boldsymbol \mu})(q)$ and a value at $q$ of an automorphism of $\widetilde 
\chi$.  Without loss of generality we may assume that $\{x_{i \alpha} 
\colon 1 \le \alpha \le r\}$ and $\{u_{j \beta} \colon 1 \le \beta \le 
r\}$ consist of the standard basis vectors for each $i$ and $j$.  Then in 
the partial fraction expansion \eqref{solution} of $T(p)$ we have that 
$\Gamma$ has the block matrix form
\begin{equation}  \label{Gammafull}
 \Gamma = -[K(\widetilde{\chi}; \lambda^i, \mu^j)]_{i=1, \dots, n_0; j=1, 
\dots, n_\infty},
\end{equation}
and that $K_{ {\boldsymbol \mu},{\bold u}}(p)$ and $K^{\bold x,
{\boldsymbol \lambda}}(q)$ 
are block row and column matrices respectively
 \begin{gather}
K_{{\boldsymbol \mu}, {\bold u}}(p) = 
\begin{bmatrix} K(\widetilde{\chi}; p, \mu^1)  & \dots & 
K(\widetilde{\chi}; p, \mu^{n_\infty})  \end{bmatrix}, \\
K^{{\bold x},{\boldsymbol \lambda}}(q) = \begin{bmatrix} 
K(\widetilde{\chi}; \lambda^1, q) \\
\vdots \\  K(\widetilde{\chi}; \lambda^{n_0}, q)\end{bmatrix}.
\end{gather}
Equating this multiplicative formula to the partial fraction expansion 
leads to the same result as in formula (2.16) in \cite{newfay}.
 
\section{Determinantal representations of algebraic curves and kernel 
bundles via Cauchy kernels} \label{S:detrep}

In \cite{hip} zero-pole interpolation problems of the sort discussed here 
were studied in a more concrete setting of vector bundles over an 
algebraic curve embedded in projective space with fiber space given as 
the kernel of a two-variable matrix pencil.  In this section we make the 
connections between that setting and the abstract compact 
Riemann surface setting 
of Section~\ref{S:absint} of this paper explicit.  As we shall see, the 
link between 
the two settings is provided by the Cauchy kernels introduced in 
Section~\ref{S:cauchyker}.

We first review the setting from \cite{hip}.  Suppose that we are given 
three $M \times M$ matrices $\sigma_1, \sigma_2, \gamma$ and let 
$U_0(z) = U_0(z_1, z_2)$ be the two-variable linear 
matrix pencil
\[
U_0(z) = z_1 \sigma_2 - z_2 \sigma_1 + \gamma, \quad
z=(z_1, z_2).
\]
We will also often consider the homogenization $U(\mu)$ (where $\mu = 
[\mu_0,\mu_1,\mu_2]$ are projective coordinates in ${\bold P}^2$) given by
\[
U(\mu) = \mu_0 U_0(\frac{\mu_1}{\mu_0}, \frac{\mu_2}{\mu_0})
=\mu_1 \sigma_2 - \mu_2 \sigma_1 + \mu_0 \gamma.
\]
Although $\det U(\mu)$ is not well-defined as a function of the 
projective variable $\mu=[\mu_0,\mu_1,\mu_2]$, nevertheless its zero set 
is well-defined and defines a curve $C \subset {\bold P}^2$ by
\[
C=\{\mu=[\mu_0,\mu_1,\mu_2] \in {\bold P}^2\colon \det U(\mu)=0\}.
\]
We shall assume that $U(\mu)$ defines a maximal irreducible determinantal 
representation of rank $r$; this means that $\det U(\mu) = F(\mu)^r$ 
where $F$ is an irreducible homogeneous polynomial of degree $m$ (so 
$M=rm$), and that $\ker U(\mu) = r$ for all smooth points $\mu$ of $C$, 
i.e., points $\mu^0$ where at least one of $\frac{\partial F}{\partial 
\mu_0}(\mu^0)$, $\frac{\partial F}{\partial \mu_1}(\mu^0)$, 
and $\frac{\partial 
F}{\partial \mu_2}(\mu^0)$ is not zero.  In case of a singular point 
$\mu^0$, we assume that $\dim \ker U(\mu^0)$ is as large as possible, 
namely $sr$ where $s$ is the multiplicity of $\mu^0$.  Under these 
conditions $E(\mu) = \ker U(\mu)$ lifts to a vector bundle $E$ of rank 
$r$ over the normalizing Riemann surface $X$ of $C$; note that the bundle 
$E$ is realized concretely as a rank $r$ subbundle of the trivial bundle 
of rank $M$ over $X$.  The normalizing Riemann surface $X$ is a 
Riemann surface 
such that there is a holomorphic mapping $\pi \colon X \to {\bold P}^2$ 
whose image equals $C$ such that $\pi$ is a one-to-one immersion on the 
inverse image of smooth points of $C$; we call $\pi \colon X \rightarrow C$ 
a birational embedding of $X$ in ${\bold P}^2$.
  For more details, see \cite{hip}.
As in \cite{hip}, we shall assume for simplicity that all the singular 
points of $C$ are nodes (i.e., $\pi^{-1}(q) = \{p^1, p^2 \}$ where $p^1$ 
and $p^2$ are distinct points on $X$ with neighborhoods $U_1$ and $U_2$ 
such that $\pi$ is an immersion at both $p^1$ and $p^2$ and the analytic 
arcs $\pi(U_1) $ and $\pi(U_2)$ intersect transversally at $q$).
We also assume that the line at infinity $\{\mu_0 = 0\}$ is nowhere 
tangent to $C$.

The holomorphic vector bundle $E_\ell$ which is dual to $E$ can be 
realized concretely as a subbundle of the trivial rank $M$ bundle over $X$ 
(with fibers now written as row vectors) via
\[
E_\ell(\mu) = \ker_\ell U(\mu)
=\{x \in {\bold C}^{1 \times M} \colon x U(\mu) = 0\}.
\]
A concrete pairing between $E_\ell \otimes {\cal O}(1) \otimes \Delta$ and 
$E \otimes {\cal O}(1) \otimes \Delta$ is given by
\begin{equation} \label{pairing}
\{u_\ell, u\} = \frac{u_\ell}{\mu_0} \frac{\xi_1\sigma_1 + \xi_2 
\sigma_2}{\xi_1\ d\lambda_1 + \xi_2\ d\lambda_2} \frac{u}{\mu_0}.
\end{equation}
Here $u_\ell$ and $u$ are local holomorphic sections of $E_\ell \otimes 
{\cal O}(1) \otimes \Delta$ and $E\otimes {\cal O}(1) \otimes \Delta$ 
respectively, $\lambda_1$ and $\lambda_2$ are meromorphic functions on 
$X$ given by $\lambda_1 = z_1 \circ \pi$ and $\lambda_2 = z_2 \circ \pi$,
and $\xi_1,\xi_2$ are arbitrary (not both zero) 
complex parameters.

If $E$ and $E_\ell$ are right and left kernel bundles determined by a 
rank $r$ maximal determinantal representation $U(\mu)$ of a curve $C$ as 
above, then it can be shown that necessarily $E \otimes {\cal O}(1) 
\otimes \Delta$ is isomorphic to a flat bundle $\chi$ over $X$ with the 
property that $h^0(\chi \otimes \Delta) = 0$, and that $E_\ell \otimes 
{\cal O}(1) \otimes \Delta$ is isomorphic to the dual $\chi^\vee$ of $\chi$.  
This isomorphism of $E \otimes {\cal O}(1) \otimes \Delta$ with $\chi$
is implemented 
explicitly by a {\it matrix of normalized sections} $u^\times(p)$.  
Explicitly, $u^\times$ is an $M \times r$ matrix whose columns are 
meromorphic sections of the pullback of $E \otimes \Delta$ to the 
universal cover $\widetilde{X}$ of $X$ such that:
\begin{enumerate}
\item[1.] $\dfrac{1}{\sqrt{dt}(R \widetilde{p})} u^\times(R \widetilde{p}) = 
\frac{1}{\sqrt{dt}(\widetilde{p})}u^\times(\widetilde{p})
 \chi^{-1}(R)$ for all $\widetilde{p} 
\in \widetilde{X}$ and all $R \in \text{Deck}(\widetilde{X}/X)$ $ \cong 
\pi_1(X)$, where $t$ is a local parameter on $X$ and $\sqrt{dt}$ is the 
corresponding local holomorphic frame for $\Delta$ lifted to the 
neighborhoods of $\widetilde{p}$ and $R\widetilde{p}$ on $\widetilde{X}$.
\item[2.] Each column of $u^\times$ has first order poles at (the points 
of $\widetilde{X}$ over) the points of $C$ at infinity, and is 
holomorphic everywhere else.
\item[3.] For each $p \in X$, the columns of $u^\times(\widetilde{p})$ 
form a basis for the fiber $(E \otimes \Delta)(p)$, where $\widetilde{p} 
\in \widetilde{X}$ is over $p$ (if $p$ is a point of $C$ at infinity we 
have first to multiply $u^\times$ by a local parameter centered at $p$).
\end{enumerate}
Simply speaking, $u^\times$ consists of a multiplicative $\Delta$-valued 
meromorphic frame for $E$, normalized to have poles exactly at the points 
of $C$ at infinity.  An isomorphism $\chi \to E \otimes {\cal O}(1) 
\otimes \Delta$ is now given explicitly by $y \to \mu_0 u^\times y$ 
where $y$ is a 
local holomorphic section of $\chi$. An $r \times M$ matrix of normalized 
sections $u_\ell^\times$ of $E_\ell$, whose rows are meromorphic sections 
of the pullback of $E_\ell \otimes \Delta$ to the universal covering 
$\widetilde{X}$ of $X$, is defined similarly, with item (1) replaced by
\begin{enumerate}
\item[$1_\ell .$] $\dfrac{1}{\sqrt{dt}}
(R \widetilde{p}) u_\ell^\times(R \widetilde p) = 
\frac{1}{\sqrt{dt}(\widetilde{p})} \chi(R) u_\ell^\times(\widetilde{p})$ 
for all $\widetilde{p} \in \widetilde{X}$ and all $R \in 
\text{Deck}(\widetilde{X}\backslash X) \cong \pi_1(X)$, where $t$ and 
$\sqrt{dt}$ are as before. 
\end{enumerate}
 An isomorphism $\chi^\vee \to E_\ell \otimes 
{\cal O}(1) \otimes \Delta$ is given explicitly by $x \to \mu_0 x^T 
u_\ell^\times$, where $x$ is a local holomorphic section of $\chi^\vee$.  
 Given $u^\times$, the {\it dual } 
matrix of normalized section $u_\ell^\times$ is determined uniquely by
\[
u_\ell^\times \frac{\xi_1 \sigma_1 + \xi_2 \sigma_2}{\xi_1\ d\lambda_1 + 
\xi_2\ d\lambda_2} u^\times = I_r
\]
(where $I_r$ is the $r \times r$ identity matrix), so that under the 
isomorphisms $\chi \cong E \otimes {\cal O}(1) \otimes \Delta$ and 
$\chi^\vee \cong E_\ell \otimes {\cal O}(1) \otimes \Delta$ the natural 
duality pairing between $\chi^\vee$ and $\chi$ equals the pairing 
\eqref{pairing}.

If we now define $K(\chi; \widetilde{p}, \widetilde{q})$ by
\begin{equation} \label{cauchypairing}
K(\chi; \widetilde{p}, \widetilde{q}) = u_\ell^\times(\widetilde{p})
\frac{\xi_1 \sigma_1 + \xi_2 \sigma_2}{\xi_1(\lambda_1(\widetilde{p}) - 
\lambda_1(\widetilde{q})) + \xi_2 (\lambda_2(\widetilde{p}) - 
\lambda_2(\widetilde{q}))} u^\times(\widetilde{q}),
\end{equation}
then $K$ has all the properties of the Cauchy kernel as defined in 
Section \ref{S:cauchyker}.  This method of constructing the Cauchy 
kernel, via a dual pair of normalized sections of the kernel bundles 
associated with a maximal determinantal representation of an algebraic 
curve $C$ embedded in ${\bold P}^2$ which has the Riemann surface $X$ as 
its normalizing surface, was presented in \cite{hip}.

Here we wish to make explicit the reverse path.  We start with a 
compact Riemann 
surface $X$ and a flat holomorphic vector bundle $\chi$ over $X$ for which 
$h^0(\chi \otimes \Delta) = 0$.  We assume as given the associated Cauchy 
kernel as developed in Section \ref{S:cauchyker}.  We then construct a 
birational embedding of $X$ into ${\bold P}^2$ with image equal to the 
curve $C$ together with a rank $r$ maximal determinantal representation 
of $C$ in such a way that we recover the Cauchy kernel $K(\chi; \cdot, 
\cdot)$ from a dual pair of normalized sections for the associated left 
and right kernel bundles associated with this determinantal 
representation of $C$, as in \eqref{cauchypairing}.  

We first need some preparations.  Let $\chi$ be a flat vector bundle over 
the Riemann surface $X$ such that $h^0(\chi \otimes \Delta)=0$.  
In addition choose two scalar meromorphic functions $\lambda_1$, 
$\lambda_2$ on $X$ such that ${\cal M}(X) = {\bold C}(\lambda_1, 
\lambda_2)$, i.e., rational functions in $\lambda_1,\lambda_2$ generate 
the whole field of (scalar) meromorphic functions on $X$.  Assume that 
all poles of $\lambda_1$ and $\lambda_2$ are simple, and denote the set 
of poles by $x^1, \dots, x^m \in X$.  Define complex numbers $c_{ik}$ ($1 
\le i \le m$, $k=1,2$) by
\[
c_{ik} = -\text{Res}_{p=x^i}\lambda_k(p)
\]
where the residue is with respect to some fixed local coordinate $t^i = 
t^i(p)$ centered at $p=x^i$.  On occasion we shall also need the next 
coefficient $-d_{ik}$ in the Laurent expansion of $\lambda_k$ at $x^\i$:
\[
\lambda_k(p) = -\frac{c_{ik}}{t^i} - d_{ik} + O(|t^i|).
\]
Define $M \times M$ matrices (where $M = mr$) $\sigma_1, \sigma_2,\gamma$ by
\begin{equation} \label{pencilcoef}
\sigma_1 = \underset{1\le i \le m}{\text{diag.}} (c_{i1} I_r), \quad
\sigma_2 = \underset{1\le i \le m}{\text{diag.}}(c_{i2}I_r) \quad
\gamma = [\gamma_{ij}]_{i,j=1,\dots,m}
\end{equation}
where
\[
\gamma_{ij} = \begin{cases} 
d_{i1}c_{i2} -d_{i2} c_{i1}, & i=j \\
(c_{i1} c_{j2}- c_{j1}c_{i2}) \dfrac{K(\chi; x^i, x^j)}{dt^j(x^j)}, & i \ne j.
\end{cases}
\]
Also define
\begin{equation} \label{normsec}
u^\times(p) = \begin{bmatrix} K(\chi; x^1, p) \\ \vdots \\ K(\chi; x^m,p) 
\end{bmatrix}, \quad
u_\ell^\times(p) = -\begin{bmatrix} K(\chi; p, x^1) & \dots & K(\chi; p, x^m) 
\end{bmatrix}
\end{equation}
Then we have the following result.

\begin{theorem} \label{T:detrep}  Let $\chi$ be a flat vector bundle over 
the Riemann surface $X$ such that $h^0(\chi \times \Delta) = 0$ with 
associated Cauchy kernel $K(\chi; \cdot, \cdot)$ and use a pair of 
meromorphic functions $\lambda_1(p),\lambda_2(p)$ on $X$  which generate 
the field ${\cal M}(X)$ of meromorphic functions on $X$ to define matrices 
$\sigma_1$, $\sigma_2$ and $\gamma$ as in \eqref{pencilcoef}.  Then:

(i)  The map $\pi_0 \colon X \to {\bold C}^2$ given by
\[
\pi_0(p) = (\lambda_1(p), \lambda_2(p))
\]
maps $X \backslash \{x^1, \dots, x^m\}$ onto the affine part $C_0$ of 
an algebraic curve $C \subset {\bold P}^2$  and extends to a 
birational embedding $\pi \colon X \rightarrow C$ of $X$ in ${\bold P}^2$. 
The defining irreducible
homogeneous polynomial $F(\mu_0, \mu_1, \mu_2)$ of $C$ is such that
 $ \det (\mu_1 \sigma_2 - 
\mu_2 \sigma_1 + \mu_0 \gamma) = F(\mu_0, \mu_1, \mu_2)^r$. 

(ii)  Denote by $E$ the kernel bundle over $C$ given in affine 
coordinates by
\[
E(\lambda) = \ker (\lambda_1 \sigma_2 - \lambda_2 \sigma_1 + \gamma),\quad
\lambda= (\lambda_1, \lambda_2).
\]
Then $\chi \cong E \otimes {\cal O}(1) \otimes \Delta$ with $u^\times$ and 
$u_\ell^\times$ given by \eqref{normsec} equal to the dual matrices of 
normalized sections of $E$ and $E_\ell$.
\end{theorem}

We shall prove Theorem \ref{T:detrep} under the assumption that all the 
singular points of $C$ are nodes.

\begin{pf}
We define the curve $C$ as the compactification in projective space of
the image
of the map $\pi$ in the statement of the theorem
\[
 C_0 = \{(\lambda_1(p), \lambda_2(p))\colon p \in X\}.
 \] 
Then $X$ is the normalizing Riemann surface for the curve $C$ and the 
degree of $C$ is equal to the number of intersections with the line at 
infinity, namely, deg $C=m$. Let $f(z_1,z_2)=0$ be an irreducible 
polynomial of degree $m$ such that $f(z_1,z_2)=0$ is the defining equation 
for $C$ (in affine coordinates).  Thus $f(\lambda_1(p),\lambda_2(p))=0$ 
for all $p \in X$.  We must show that $z_1\sigma_2 - z_2 \sigma_1 + 
\gamma$ is a maximal determinantal representation of $f(z_1,z_2)^r=0$, 
that $E \otimes {\cal O}(1) \otimes \Delta \cong \chi$ and that 
$u^\times$ and $u_\ell^\times$ are dual matrices of normalized sections 
of $E$ and $E_\ell$.

The first step is to prove the identities
\begin{align}
(\lambda_1(p) \sigma_2 -\lambda_2(p)\sigma_1+\gamma)u^\times(p) & =0 
\label{identity1} \\
u_\ell^\times(p) (\lambda_1(p) \sigma_2 - \lambda_2(p) \sigma_1 + \gamma) 
& = 0 
\label{identity2} \\
\frac{u_\ell^\times(p) (\xi_1 \sigma_1 + \xi_2 \sigma_2 ) u^\times(p)}
{\xi_1 \ d\lambda_1(p) + \xi_2\ d\lambda_2(p)} & = 1.
\label{identity3}
\end{align}

To prove \eqref{identity1}, set 
\[
h(p) = (\lambda_1(p) \sigma_2 - \lambda_2(p) \sigma_1 + \gamma) 
u^\times(p).
\]
Note that $h(p)$ is a meromorphic section of $\chi \otimes \Delta$.  To 
check that $h=0$ it suffices to check that $h$ has no poles (since 
$h^0(\chi \otimes \Delta)=0$).  The only possible poles in the formula 
for $h$ occur at the points $x^1, \dots, x^m$ and these are at most 
double poles.  For $\alpha = 1, \dots, m$, let us write down the Laurent 
expansion for $h$ near $x^\alpha$ as
\[
h(p) = [h]^{\alpha,-2}(t^\alpha)^{-2} + [h]^{\alpha, -1} (t^\alpha)^{-1} 
+ [\text{analytic at $x^\alpha$]}.
\]
We must show that $[h]^{\alpha, -2}=0$ and $[h]^{\alpha, -1} = 0$ 
for $1 \le 
\alpha \le m$.

Each of $[h]^{\alpha, -2}$ and $[h]^{\alpha, -1}$ in turn is a block $m 
\times 1$ column matrix: $[h]^{\alpha, -2} = \left[ [h]^{\alpha, -2}_i\right]$ 
and
$[h]^{\alpha, -1} = \left[ [h]^{\alpha, -1}_i \right]$ with $i = 1, 
\dots m$.   We compute
\begin{align} \notag
[h]^{\alpha, -2}_i =& \sum_{j=1}^m (-c_{\alpha 1} c_{i2} \delta_{ij} + 
c_{\alpha 2} c_{i 1} \delta_{ij}) \delta_{j \alpha}\cdot 
(- dt^\alpha(x^\alpha))\\
\notag
=& ( c_{\alpha 1} c_{\alpha 2} - c_{\alpha 2} c_{\alpha 1}) 
dt^\alpha(x^\alpha) = 0
\end{align}
where $\delta_{ij}$ is the Kronecker delta.
Similarly,
\begin{align} \notag
[h]^{\alpha, -1}_i =& \sum_{j,\ j\ne \alpha} (-c_{\alpha 1} c_{i 2} 
\delta_{ij} +c_{\alpha 2} c_{i 1} \delta_{ij}) K(\chi, x^j, x^\alpha)  \\
\label{hresidue}
 & + \sum_j^m \{-d_{\alpha 1} c_{i 2} \delta_{ij} + d_{\alpha 2} c_{i 1} 
\delta_{ij} + \gamma_{ij} \} \delta_{j \alpha}(- dt^\alpha(x^\alpha)).
\end{align}
For $i = \alpha$ \eqref{hresidue} becomes
\[
(-d_{\alpha 1} c_{\alpha 2} + d_{\alpha 2} c_{\alpha 1} + d_{\alpha 1} 
c_{\alpha 2} - c_{\alpha 1} d_{\alpha 2})(- dt^\alpha(x^\alpha))= 0
\]
while, for $i \ne \alpha$, \eqref{hresidue} becomes
\[
(-c_{\alpha 1} c_{i2} + c_{\alpha 2} c_{i 1}) K(\chi; x^i, x^\alpha) + 
(c_{i1} c_{\alpha 2} - c_{\alpha 1} c_{i 2}) \frac{K(\chi; x^i, 
x^\alpha)}{dt^\alpha(x^\alpha)}(- dt^\alpha(x^\alpha)) = 0.
\]
Thus $h=0$ and \eqref{identity1} follows.

By a similar calculation of Laurent series coefficients one can verify 
\eqref{identity2}. 

To verify \eqref{identity3},
by a standard lemma (see Proposition 2.3 in \cite{hip}), it suffices to 
show that
\[
\frac{u^\times_\ell(p) \sigma_k u^\times(p)}{d\lambda_k(p)} = 1 \text{ 
for } k=1,2.
\]
The numerator of the expression on the left is given by
\begin{gather} \notag
\begin{bmatrix} K(\chi; p,x^1) & \dots & K(\chi; p,x^m) \end{bmatrix}
\begin{bmatrix} -c_{ik} & & \\ & \ddots & \\  & & - c_{mk} \end{bmatrix}
\begin{bmatrix} K(\chi; x^1,p) \\ \vdots \\ K(\chi; x^m, p) 
\end{bmatrix}  \\
\notag
= -\sum_{i=1}^m K(\chi; p, x^i) [\text{Res}_{p=x^i} \lambda_k(p)] K(\chi; 
x^i,p)
\end{gather}
and hence
\[
\frac{u^\times_\ell(p) \sigma_k u^\times(p)}{d\lambda_k(p)} = 
- \frac{\sum_{i=1}^m K(\chi; p, x^i) [\text{Res}_{p=x^i}\lambda_k(p)] 
K(\chi; x^i,p)}
{d\lambda_k(p)}.
\]
The double pole at each $x^i$ in the numerator is cancelled by a double 
pole at each $x^i$ in the denominator.  Note also that the product of two 
half-order differentials in the numerator is cancelled by the 
differential in the denominator.  The resulting quotient is a 
well-defined holomorphic section of $Hom(\chi,\chi)$ which has the value 
$I$ at $x^1, \dots, x^m$, and hence must equal $I$ at all $p \in X$.  
Equation \eqref{identity3} now follows.

Denote by $d(z_1,z_2)$ the polynomial $d(z_1,z_2)=\det 
(z_1\sigma_2-z_2\sigma_1 + \gamma)$.  By \eqref{identity1} or 
\eqref{identity2} we know that $d(\lambda_1(p),\lambda_2(p))=0$ 
identically in $p \in X$. Since $f$ is by assumption the irreducible 
defining polynomial for $X$, it follows that $f(z_1,z_2) | d(z_1,z_2)$, 
and hence,
\begin{equation} \label{fact}
d(z_1,z_2) = f(z_1,z_2)^s g(z_1,z_2)
\end{equation}
for some positive integer $s$ and some polynomial $g$ relatively prime 
with respect to $f$.  By inspection we see that $d$ has 
degree equal to $M=mr$, 
while, as already mentioned, $f$ has degree equal to $m$.  From the 
factorization we see that deg $d \ge s(\text{deg } f)$, or
$r \ge s$.  On the other hand, at a smooth point $(z_1,z_2) \in C$, we 
know that $\dim E(z) \ge r$ since the $r$ linearly independent columns of 
$u^\times(p)$ are in $\ker (\lambda_1(p) \sigma_2 - \lambda_2(p) \sigma_1 
+ \gamma)$.  But also, from the factorization \eqref{fact} and an 
inductive argument working with the minors the matrix pencil 
$z_1 \sigma_2 - z_2 \sigma_1 + \gamma$ (see the proof of 
Theorem 3.2 in \cite{hip}) one can show that $\dim E(z) \le s$.
From $r \ge s$ and $r \le \dim E(z) \le s$ we conclude that $r=s$.
From \eqref{fact} and degree counting we conclude that the polynomial $g$
is a constant, and without loss of generality, $d=f^r$.  Thus 
$z_1\sigma_2 - z_2 \sigma_1 + \gamma$ is a maximal determinantal 
representation of $f(z_1,z_2)^r=0$ as asserted, except that we still must 
 check that for any node $q$ on $C$ the columns of $u^\times(p^1)$ and 
$u^\times(p^2)$ are linearly independent, where $\pi^{-1}(q) = \{p^1, p^2\}$.
It then follows from \eqref{identity1}, \eqref{identity2} and 
\eqref{identity3} that $u^\times(p)$ and $u^\times_\ell(p)$ form the 
associated dual pair of normalized cross-sections for $E$ and $E_\ell$ 
respectively.

We first  claim if $L$ is a straight line nowhere tangent to $C$ and 
$y^1, \dots, y^m$ are the preimages on $X$ of the points of intersection 
of $C$ with $L$, then the block matrix $[K(\chi; x^i, y^j]_{i,j = 1, 
\dots, m}$ is invertible.   This follows immediately from the discussion 
of the full rank zero-pole interpolation problem at the end of Section 4,
as we now show. 
Since the divisor $x^1 + \dots + x^m - y^1 - \dots - y^m$ is equivalent 
to $0$, the input bundle in the full rank zero-pole 
interpolation problem with zeros $x^1, \dots, x^m$ and poles $y^1, \dots, y^m$
and output bundle $\chi$ is again (isomorphic to) $\chi$. Since $h^0(\chi 
\otimes \Delta) = 0$, the matrix $\Gamma$ \eqref{Gammafull} is invertible.
Next, by taking any line $L$ through the node $q$ which is nowhere tangent to 
$C$, we see that the columns of $u^\times(p^1)$ and $u^\times(p^2)$ are 
columns in a $M \times M$ invertible matrix (namely, the associated 
matrix $\Gamma$), and hence are linearly 
independent. 
\end{pf}

We remark that the same proof works when the singularities of $C$ are 
any ordinary singular points, or more generally, are such that a singular 
point of multiplicity $s$ has $s$ distint preimages on $X$. 

\section{The concrete interpolation problem for meromorphic bundle maps 
between kernel bundles of determinantal representations 
of an algebraic curve} \label{S:conint}

In the paper \cite{hip} the following problem was considered.  We are 
given an irreducible algebraic curve $C$ in ${\bold P}^2$ together with 
its normalizing compact Riemann surface $X$ and the normalization map 
$\pi \colon X \to C$.  We assume that the defining polynomial for $C$ is 
an irreducible polynomial $f$ of degree $m$ (in affine coordinates).  For 
simplicity we assume again that the only singularities of $C$ are nodes and 
that $C$ intersects the line at infinity in $m$ distinct smooth points.  
We suppose in addition that $f^r$ has a maximal determinantal representation
\[
f^r(z_1,z_2) = \det (z_1 \sigma_2 - z_2 \sigma_1 + \widetilde{\gamma})
\]
where $\sigma_1$, $\sigma_2$ and $\widetilde{\gamma}$ are $M \times M$ 
matrices ($M=mr$), with which is associated the kernel bundle 
$\widetilde{E}$ of rank $r$ over $C \backslash C_{sing}$ ($C_{sing}$ is 
the set of the singular points of $C$) with fibers (over affine points) 
given by 
\begin{equation}  \label{outputbundle}
\widetilde{E}(z) = \ker (z_1 \sigma_2 - z_2 \sigma_1 + 
\widetilde{\gamma}).
\end{equation}
As explained in Section \ref{S:detrep}, we may consider the pullback of 
$\widetilde{E}$ to $X \backslash \pi^{-1}(C_{sing})$ as extended to a 
rank $r$ vector bundle over all of $X$.  We also have the left kernel 
bundle $\widetilde{E}_\ell$ where
\begin{equation}  \label{dualoutputbundle}
\widetilde{E}_\ell(z) = \ker_\ell (z_1 \sigma_2 - z_2 \sigma_1 + 
\widetilde{\gamma})
\end{equation}
with pullback under $\pi$ also extendable to a rank $r$ vector bundle 
defined over all of $X$.  These bundles, or more precisely their twists 
$\widetilde{E} \otimes {\cal O}(1) \otimes \Delta$ and 
$\widetilde{E}_\ell \otimes {\cal O}(1) \otimes \Delta$, have the canonical 
pairing \eqref{pairing} with each other, as explained in Section 
\ref{S:detrep}.

The data for the concrete interpolation (CONINT) 
problem to be considered in this 
section consists of:
\begin{enumerate}
\item[(D1)] $n_\infty$ distinct smooth, finite points $\mu^1 = (\mu^1_1, 
\mu^1_2), \dots, \mu^{n_\infty} = (\mu^{n_\infty}_1, \mu^{n_\infty}_2)$ 
of $C$ (the preassigned poles),
\item[(D2)] for each $j=1, \dots, n_\infty$, a linearly independent set 
$\{\varphi_{j1}, \dots, \varphi_{j,s_j} \}$ of $s_j$ vectors in the fiber 
$\widetilde{E}(\mu^j)$ (the preassigned pole vectors),
\item[(D3)] $n_0$ distinct smooth, finite points $\lambda^1 = (\lambda^1_1, 
\lambda^1_2), \dots, \lambda^{n_0}=(\lambda^{n_0}_1, \lambda^{n_0}_2)$ of 
$C$ (the preassigned zeros),
\item[(D4)] for each $i=1, \dots, n_0$, a linearly independent set 
$\{\psi_{i1}, \dots, \psi_{i,t_i}\}$ of $t_i$ vectors in the fiber 
$\widetilde{E}_\ell(\lambda^i)$ (the preassigned null vectors), and
\item[(D5)] for each pair of indices $(i,j)$ for which $\lambda^i = \mu^j =: 
\xi^{ij}$, a choice of a local coordinate $t^{ij}$ on $X$ centered at 
$\xi^{ij}$ and a collection of numbers $\{\rho_{ij, \alpha \beta} \colon 
1 \le \alpha \le t_i, 1 \le \beta \le s_j\}$ (the preassigned coupling 
numbers with respect to the chosen local coordinate).
\end{enumerate}
The interpolation problem then is to find an $M \times M$ matrix $\gamma$ 
defining a maximal determinantal representation of $f^r$
\begin{equation} \label{detrepforfr}
f^r(z_1,z_2) = \det (z_1 \sigma_2 - z_2 \sigma_1 + \gamma)
\end{equation}
giving the kernel bundle $E$ over $X$ with fiber over a smooth finite 
point $z \in C$ given by
\begin{equation} \label{inputbundle}
E(z) = \ker (z_1 \sigma_2 - z_2 \sigma_1 + \gamma)
\end{equation}
and the left kernel bundle $E_\ell$ given by
\begin{equation} \label{dualinputbundle}
E_\ell(z) = \ker_\ell (z_1 \sigma_2 - z_2 \sigma_1 + \gamma)
\end{equation}
together with meromorphic bundle maps
\[
S \colon E \to \widetilde{E},\quad S_\ell \colon \widetilde{E}_\ell \to 
E_\ell
\]
(where we write bundle maps on left kernel bundles as acting from the 
right), where $S \otimes I_{{\cal O}(1) \otimes \Delta}$ and $S_\ell 
\otimes I_{{\cal O}(1) \otimes \Delta}$ are transposes of each other with 
respect to the pairing \eqref{pairing}, so that $S$ (and $S_\ell$) act as 
the identity operator $I$ on the corresponding fibers at the points at 
infinity, and the following set of interpolation conditions is satisfied:
\begin{enumerate}
\item[(I1)] $S$ has poles only at $\mu^1, \dots, \mu^{n_\infty}$; for each 
$j=1, \dots, n_\infty$, the pole of $S$ at $\mu^j$ is simple, and the 
vectors $\{\varphi_{j1}, \dots, \varphi_{j,n_\infty}\}$ span the image 
space of the residue $R_j \colon E(\mu^j) \to \widetilde{E}(\mu^j)$ of 
$S$ at $\mu^j$.
\item[(I2)] The bundle map $S_\ell^{-1} \colon
 E_\ell \to \widetilde{E}_\ell$ has 
poles only at the points $\{\lambda^1, \dots, \lambda^{n_0}\}$; for each 
$i=1, \dots, n_0$, the pole of $S_\ell^{-1}$ at $\lambda^i$ is simple and 
the vectors $\{\psi_{i1} \dots, \psi_{i t_i}\}$ span the image space of 
the residue $\widehat{R}_i \colon E_\ell(\lambda^i) \to 
\widetilde{E}_\ell(\lambda^i)$ of $S_\ell^{-1}$ at $\lambda^i$.
\item[(I3)]  For each pair of indices $(i,j)$ where $\lambda^i = \mu^j =: 
\xi^{ij}$, and for $\alpha = 1, \dots, t_i$, let $\psi_{i \alpha}(p)$ be a 
local holomorphic section of $\widetilde{E}_\ell$ near $\xi^{ij}$ with
\[
\psi_{i \alpha}(t^{ij}) = \psi_{i \alpha} + \psi_{i \alpha 1} t^{ij} + 
o(t^{ij})
\]
such that $\psi_{i \alpha} S_\ell(p)$ has analytic continuation to $p = 
\xi^{ij}$ with value there equal to $0$.  Then, for any choice of complex 
parameters $\xi_1$ and $\xi_2$
\[
\psi_{i \alpha 1} 
\frac{\xi_1 \sigma_1 + \xi_2 \sigma_2}
{\xi_1 \lambda_1^{\prime}(\xi^{ij}) + \xi_2 
\lambda_2^{\prime}(\xi^{ij})}\varphi_{j \beta}
 - \psi_{i \alpha} (\xi_1 \sigma_1 + \xi_2 \sigma_2) \varphi_{j \beta}
\frac{\xi_1 \lambda_1^{\prime \prime}(\xi^{ij}) + \xi_2 \lambda_2^{\prime 
\prime}(\xi^{ij})}
{2(\xi_1 \lambda_1^{\prime}(\xi^{ij}) + \xi_2 
\lambda_2^{\prime}(\xi^{ij}))^2} = \rho_{ij, \alpha \beta}.
\]
Here ${}^{\prime} = \dfrac{d}{dt^{ij}}$.  
\end{enumerate}
It can be shown that a necessary consistency condition 
on the data set for the problem to have a 
solution is that 
\begin{equation} \label{ZP}
\psi_{i \alpha} (\xi_1 \sigma_1 + \xi_2 \sigma_2) \varphi_{j \beta} = 0.
\end{equation}

Given the data set (D1)--(D5) we form the  $n_0 \times n_\infty$
block matrix $\Gamma^0 = 
[\Gamma^0_{ij}]$ ($1 \le i \le n_0, 1 \le j \le n_\infty$) where 
$\Gamma^0_{ij}=[\Gamma^0_{ij,\alpha \beta}]$ ($1 \le \alpha \le t_i$, $1 
\le \beta \le s_j$) in turn is the $t_i \times s_j$ matrix with entries 
given by
\begin{equation}  \label{Gamma0} 
\Gamma^0_{ij, \alpha \beta} = \begin{cases}
\psi_{i \alpha} \frac{\xi_1 \sigma_1 + \xi_2 \sigma_2}
{\xi_1(\mu^j_1 - \lambda^i_1) + \xi_2(\mu^j_2 - \lambda^i_2)} 
\varphi_{j \beta} & \text{if } \lambda^i \ne \mu^j  \\
-\rho_{ij, \alpha \beta} & \text{if } \lambda^i = \mu^j. 
\end{cases}
\end{equation}
Additional matrices which we shall need are
\begin{align}
    A_1 = \begin{bmatrix} \mu^1_1I_{s_1} & & \\
    & \ddots & \\
    & & \mu^{n_\infty}_1 I_{s_{n_\infty}} \end{bmatrix},
& \quad
    A_2 = \begin{bmatrix} \mu^1_2I_{s_1} & & \\
    & \ddots & \\
    & & \mu^{n_\infty}_2 I_{s_{n_\infty}} \end{bmatrix},  \notag
\\
    Z_1=\begin{bmatrix}\lambda^1_1 I_{t_1} & &  \\
    & \ddots &  \\
    &        &  \lambda^{n_0}_1 I_{t_{n_0}}
    \end{bmatrix}, 
& \quad 
    Z_2=\begin{bmatrix}\lambda^1_2 I_{t_1} & &  \\
    & \ddots &  \\
    &        &  \lambda^{n_0}_2 I_{t_{n_0}}
    \end{bmatrix},    \notag
\\
    \varphi_j = \begin{bmatrix} \varphi_{j1} & \dots & \varphi_{j s_j} 
    \end{bmatrix},  
& \quad
    \varphi = \begin{bmatrix} \varphi_1 & \dots & \varphi_{n_\infty} 
    \end{bmatrix},   \notag
\\
  \psi_i = \begin{bmatrix} \psi_{i1} \\ \vdots \\ \psi_{i t_i} 
  \end{bmatrix}, 
& \quad  
  \psi = \begin{bmatrix} \psi_1 \\ \vdots \\ \psi_{n_0} \end{bmatrix}.
  \label{matrices}
\end{align}

The solution of the concrete interpolation problem (CONINT) obtained in 
\cite{hip} is as follows.

\begin{theorem}  \label{T:conint} (See Theorem 4.1 of \cite{hip}.)
Assume that we are given a curve $C$ 
with defining irreducible polynomial $f$, a maximal determinantal 
representation for $f^r$ as in \eqref{detrepforfr} together with 
associated kernel bundle $\widetilde{E}$ \eqref{outputbundle}
 and left kernel bundle 
$\widetilde{E}_\ell$ \eqref{dualoutputbundle}
and a data set (D1)--(D5) for the interpolation 
problem (I1)--(I3).  Then the interpolation problem has a solution if and 
only if the interpolation data satisfy the compatibility conditions 
\eqref{ZP} 
at the overlapping zeros and poles, and the matrix $\Gamma^0$ given by 
\eqref{Gamma0} is square and invertible.  
In this case the unique solution of the 
interpolation problem (I1)--(I3) is given by
\begin{equation}  \label{gamma}
\gamma = \widetilde{\gamma} - \sigma_1 \varphi (\Gamma^0)^{-1} \psi \sigma_2 
+ \sigma_2 \varphi (\Gamma^0)^{-1} \psi \sigma_1
\end{equation}
with associated kernel bundle $E$ \eqref{inputbundle} and left kernel 
bundle $E_\ell$ \eqref{dualinputbundle},
 with $S(z)$ given by
\begin{equation} \label{formulaforS}
S(z) = [I + \varphi (\xi_1(z_1I - A_1) + \xi_2(z_2 I -A_2))^{-1} 
(\Gamma^0)^{-1} \psi (\xi_1 \sigma_1 + \xi_2 \sigma_2)]|_{E(z)}
\end{equation}
and with $S_\ell^{-1}(z)$ given (as a right multiplication operator) by
\begin{equation} \label{formulaforSell}
S_\ell^{-1}(z) = [I-(\xi_1 \sigma_1 + \xi_2 \sigma_2) \varphi 
(\Gamma^0)^{-1} \left(\xi_1 (z_1 I - Z_1) + \xi_2( z_2 I - Z_2) \right)^{-1}
\psi] |_{E_\ell(z)}.
\end{equation}
Here the matrices $A_1,A_2,Z_1,Z_2,\psi,\varphi$ are as in 
\eqref{matrices}.
\end{theorem}

The main goal of this section is to use the machinery developed in 
Section \ref{S:detrep} to make explicit the connections between Theorem 
\ref{T:conint} and Theorem \ref{T:absint} of Section \ref{S:absint}.
  
Suppose therefore that $\chi$ and $\widetilde{\chi}$ are two flat bundles 
over the Riemann surface $X$ with $h^0(\chi \otimes 
\Delta)=0=h^0(\widetilde{\chi} \otimes  \Delta)$.  We use a fixed pair 
$(\lambda_1(p), \lambda_2(p))$ of meromorphic functions on $X$ generating 
${\cal M}(X)$ to produce a map $\pi \colon X \to C$.  The respective 
Cauchy kernels $K(\chi; \cdot, \cdot)$ and $K(\widetilde{\chi}; \cdot, 
\cdot)$  generate corresponding maximal determinantal representations 
$z_1 \widetilde{\sigma}_2 - z_2 \widetilde{\sigma}_1 + 
\widetilde{\gamma}$ and $z_1 \sigma_2^{\prime} - z_2 \sigma_1^{\prime} + 
\gamma^{\prime}$ for $f(z_1,z_2)^r=0$ as in Theorem \ref{T:detrep}.  Note 
however that the formulas for $\sigma_1$ and $\sigma_2$ in 
\eqref{pencilcoef} depend only on the choice of embedding functions 
$\lambda_1(p)$ and $\lambda_2(p)$, and not on the particular flat bundle 
$\chi$; hence we may and shall write simply $\sigma_i$ in place of 
$\widetilde{\sigma}_i$ and $\sigma_i^{\prime}$ for $i=1,2$.  We also have 
associated dual pairs of matrices of normalized sections:  
$\widetilde{u}^\times, \widetilde{u}_\ell^\times$ for $\widetilde{E}=
\ker(z_1 \sigma_2 - z_2 \sigma_1 + \widetilde{\gamma})$ and 
$\widetilde{E}_\ell = \ker_\ell (z_1 \sigma_2 - z_2 \sigma_1 + 
\widetilde{\gamma})$ respectively, and $u^{\times \prime}, u^{\times 
\prime}_\ell$ for $E= \ker (z_1 \sigma_2 - z_2 \sigma_1 + 
\gamma^{\prime})$ and $E_\ell = \ker_\ell (z_1 \sigma_2 - z_2 \sigma_1 + 
\gamma^{\prime})$.  Since $\widetilde{u}^\times$ implements an 
isomorphism between $\widetilde{\chi}$ and $\widetilde{E} \otimes {\cal 
O}(1) \otimes \Delta$ and $u^{\times \prime}$ implements an isomorphism 
between $\chi$ and $E^\prime \otimes {\cal O}(1) \otimes \Delta$, any 
meromorphic bundle map $T \colon \chi \to \widetilde{\chi}$ induces a 
meromorphic bundle map $S \colon E^\prime \to \widetilde{E}$ determined by
\[
S(p) u^{\times \prime}(p) = \widetilde{u}^\times(p) T(p).
\]
However, in the solution of (CONINT) from \cite{hip} stated in Theorem 
\ref{T:conint}, the solution $S$ is normalized to act as the identity 
operator over the points of $C$ at infinity.  In order for the map $S$ 
constructed as above from the abstract bundle map $T \colon \chi \to 
\widetilde{\chi}$ to achieve this normalization, we must make an adjustment
\begin{equation} \label{adjustment}
\alpha (z_1 \sigma_2 - z_2 \sigma_1 + \gamma^\prime) \beta = z_1 \sigma_2 
- z_2 \sigma_1 + \gamma
\end{equation}
on the input determinantal representation, where $\alpha, \beta \in 
GL(M^{mr}, {\bold C})$.  If $u^\times, u^\times_\ell$ is the dual pair of 
normalized sections for $E=\ker (z_1 \sigma_2 - z_2 \sigma_1 + \gamma)$ 
and $E_\ell = \ker_\ell (z_1 \sigma_2 - z_2 \sigma_1 + \gamma)$, then
\[
u^\times(p) = \beta^{-1} u^{\times \prime}(p), \quad u^\times_\ell(p) = 
u^{\times \prime}_\ell(p) \alpha^{-1}
\]
and we seek to solve instead
the equation
\[
S(p) u^\times(p) = \widetilde{u}^\times(p) T(p),
\]
or equivalently
\begin{equation} \label{intertwining}
S(p) \beta^{-1} u^{\times \prime}(p) = \widetilde{u}^\times(p) T(p)
\end{equation}
for $S,\alpha,\beta$ subject to the proviso that $S(x^i)= 
I_{E(x^i)=\widetilde{E}(x^i)}$ for $i=1, \dots, m$.  Since the columns of 
$u^{\times \prime}(p)$ and $\widetilde{u}^\times(p)$ (after multiplication 
by a local parameter on $X$ at $p$) simply form a 
standard basis in ${\bold C}^M$ when evaluated at $x^1, \dots, x^m$, we 
see that we should take
\[
\beta = \begin{bmatrix} T(x^1) & & \\ & \ddots & \\
 & & T(x^m) \end{bmatrix}^{-1}.
 \]
 In order to guarantee $\alpha \sigma_k \beta = \sigma_k$ for $k=1,2$ as 
 required in \eqref{adjustment} we then take
 \[
 \alpha = \begin{bmatrix} T(x^1) & & \\ & \ddots & \\  & & T(x^m) 
 \end{bmatrix}.
 \]
 Thus $\gamma = \alpha \gamma^{\prime} \beta$ is given by $\gamma = 
 [\gamma_{ij}]_{i,j=1,\dots, m}$ with
 \begin{equation} \label{inputgamma}
 \gamma_{ij} = \begin{cases} d_{i1}c_{i2} - c_{i1} d_{i2} & \text{if } i=j \\
 (c_{i1}c_{j2} - c_{j1}c_{i2})T(x^i)K(\chi;x^i,x^j)T(x^j)^{-1} & \text{if 
 } i \ne j.
 \end{cases}
 \end{equation}
 We remark that the ``adjustment'' of $\gamma$ by the values of a bundle 
 map at the points of $X$ over the points of $C$ at infinity plays a 
 central role in the construction of triangular models for commuting 
 nonselfadjoint operators; see \cite{vin1} and \cite{lkmv}, Chapter 12.

 We now suppose that we are given an abstract interpolation data set 
 $\omega$ as in \eqref{dataset} for an Abstract Interpolation Problem 
 (ABSINT) as in Section \ref{S:absint} with output bundle 
 $\widetilde{\chi}$, and let $y \to \widetilde{u}^\times y$ and $x \to 
 x^T \widetilde{u}^\times_\ell$ be the associated bundle isomorphisms 
 from $\widetilde{\chi}$ to $\widetilde{E} \otimes {\cal O}(1) \otimes 
 \Delta$ and from $\widetilde{\chi}^\vee$ to $\widetilde{E}_\ell \otimes 
 {\cal O}(1) \otimes \Delta$, where $\widetilde{E}$ and 
 $\widetilde{E}_\ell$ are the right and left kernel bundles respectively 
 associated with the maximal determinantal representation $z_1 \sigma_2 - 
 z_2 \sigma_1 + \widetilde{\gamma}$ (with 
 $\sigma_1,\sigma_2,\widetilde{\gamma}$ given by \eqref{pencilcoef}, 
 and with 
 $\widetilde{u}^\times$ and $\widetilde{u}^\times_\ell$ given by 
 \eqref{normsec}, all with $\widetilde{\chi}$ in place of $\chi$).  We 
 assume that the pair of meromorphic functions $\lambda_1(p)$ and 
 $\lambda_2(p)$ is chosen in such a way that the set of poles $x^1, 
 \dots, x^m$ is disjoint from the preassigned poles $\mu^1, \dots, 
 \mu^{n_\infty}$ and the set of preassigned zeros $\lambda^1, \dots, 
 \lambda^{n_0}$. Define a data set $\omega_0$ for a (CONINT) problem as 
 follows:
 \begin{enumerate}
 \item The preassigned poles consist of the points 
 $\pi(\mu^1) = (\mu^1_1, \mu^1_2),
 \dots, \pi(\mu^{n_\infty}) = (\mu^{n_\infty}_1, \mu^{n_\infty}_2)$ with 
 associated pole vectors $\varphi_{j \beta}  \in 
 \widetilde{E}(\pi(\mu^j))$ given by $\varphi_{j \beta} = 
 \widetilde{u}^\times(\mu^j) u_{j \beta}$ for $j=1, \dots, n_\infty$ and 
 $\beta = 1, \dots, s_j$.
 \item The preassigned zeros consist of the points $\pi(\lambda^1) = 
 (\lambda^1_1, \lambda^1_2), \dots, \pi(\lambda^{n_0}) = ( \lambda^{n_0}_1, 
 \lambda^{n_0}_2)$ with associated sets of null vectors $\psi_{i \alpha}
 \in \widetilde{E}_\ell(\pi(\lambda))$
 given by $\psi_{i \alpha} = x^T_{i \alpha} 
 \widetilde{u}^\times_\ell(\lambda^i)$ for $i=1, \dots, n_0$ and 
 $\alpha=1, \dots, t_i$.
 \item For those pairs of indices $(i,j)$ for which $z^i=w^j=:\xi^{ij}$
 we take the associated coupling numbers $\rho_{ij, \alpha \beta}$ to be 
 the same as those specified for the (ABSINT) problem.
 \end{enumerate}
 With this choice of data set, the reader can check that the matrices 
 $\varphi$ and $\psi$ as defined in \eqref{matrices} reduce to
 \[
 \varphi = \begin{bmatrix} K_{\boldsymbol \mu, {\bold u}}(x^1) \\ \vdots \\
 K_{\boldsymbol \mu,{\bold u}}(x^m) \end{bmatrix}, \quad
 \psi = -\begin{bmatrix} K^{{\bold x}, \boldsymbol \lambda}(x^1) & 
 \dots & K^{{\bold 
 x}, \boldsymbol \lambda}(x^m) \end{bmatrix}
 \]
 where the notation $K_{\boldsymbol \mu,{\bold u}}(p)$ and 
 $K^{{\bold x}, \boldsymbol \lambda}(p)$ 
 is as in the statement of Theorem \ref{T:absint}.
 
 It turns out that the (ABSINT) problem with data set $\omega$ is 
 equivalent to the (CONINT) problem with data set $\omega_0$ under the 
 identifications $\pi \colon X \to C$ and $\widetilde{u}^\times: 
 \widetilde{\chi} \to \widetilde{E} \otimes {\cal O}(1) \otimes \Delta$
 and $\widetilde{u}^\times_\ell: \widetilde{\chi}^\vee \to 
 \widetilde{E}_\ell \otimes {\cal O}(1) \otimes \Delta$ sketched above, 
 in that a solution $T$ of (ABSINT) corresponds to a solution $S$ of 
 (CONINT) under the correspondence (including the normalization
 at the points over infinity) between 
 abstract bundle maps $T$ and concrete bundle maps $S$ discussed above.
 (For the first two interpolation conditions, this observation is rather 
 transparent. For the third interpolation condition (I3), this requires 
 the relation between the interpolation condition (I3) for the (CONINT) 
 problem with a flat connection on the bundle $\widetilde{E}_\ell \otimes
 {\cal O}(1)  \otimes \Delta$ and the correspondence of this 
 connection with the coefficient $A_\ell(p)$ in the Laurent expansion of 
 the Cauchy kernel $K(\widetilde{\chi}; \cdot, \cdot)$ along the 
 diagonal; this is explained in Section 3.2 of \cite{hip}.)  
 Hence the formula for the solution $T$ of (ABSINT) in \eqref{solution} must 
 correspond to the formula for the solution $S$ of (CONINT) in 
 \eqref{formulaforS} under the correspondence \eqref{intertwining}.  The 
 point of the next result is to verify this directly; in addition we see 
 that the matrix $\Gamma$ appearing in Theorem \ref{T:absint} is 
 identical to the matrix $\Gamma^0$ appearing in Theorem \ref{T:conint}.
 
 \begin{theorem} \label{T:comparison}
 Let $\omega$ be the data set for an (ABSINT) problem with 
 $\omega_0$ the corresponding data set for a (CONINT) problem.  Then 
 $\Gamma = \Gamma_0$.  Furthermore, if $\Gamma$ is invertible and $\chi$ 
 is the input bundle for which (ABSINT) has a solution, then a solution 
 $T$ of (ABSINT) is related to the unique solution $S$ of (CONINT) having 
 value identity on the fibers over the points at infinity according to 
 the intertwining condition \eqref{intertwining}.
 \end{theorem}
 
 \begin{pf}  The fact that $\Gamma = \Gamma_0$ is a simple consequence of 
 the definitions and of the formula \eqref{cauchypairing} expressing the 
 Cauchy kernel $K(\widetilde{\chi}; \cdot, \cdot)$ in  terms of a dual 
 pair $\widetilde{u}^\times, \widetilde{u}^\times_\ell$ of normalized sections
 of $\widetilde{E}$ and $\widetilde{E}_\ell$.
 
 It remains to verify the intertwining relation \eqref{intertwining}
 \[
 S(p) \beta^{-1} u^{\times \prime}(p) = \widetilde{u}^\times(p) T(p)
 \]
 where $S$ is given by \eqref{formulaforS}, $T$ by \eqref{solution} and
 $\beta^{-1} = \underset{1\le i\le m}{\text{diag.}} \{T(x^i)\}$.
 We compute
 \begin{align} \notag
 [S(p) \beta^{-1} u^{\times \prime}(p)]_i = &
 \sum_{j=1}^m S_{ij}(p) T(x^j) K(\chi; x^j,p)  \\
 \notag
 =  T(x^i) K(\chi; x^i,p) & -
 K_{\boldsymbol \mu,{\bold u}}(x^i) \cdot
 \underset{i'}{\text{diag.}}
 \{( \xi_1 \lambda_1(p)+\xi_2 \lambda_2(p) -  \xi_1 \mu^{i^\prime}_1
  -\mu^{i^\prime}_2)^{-1} I_{s_{i^\prime}} \} \cdot  
  \\
  & \cdot 
  \Gamma^{-1} \cdot \sum_{j=1}^m (\xi_1 c_{j1} + \xi_2 c_{j2})
   K^{{\bold x}, \boldsymbol \lambda}( \lambda^j) 
  T(x^j)
   K(\chi; x^j,p).
 \label{Sside1}
 \end{align}
 From \eqref{solution} with $x^i$ in place of $p$ and $p$ in place of $q$ 
 (and hence $T(p)$ in place of $Q$) we see that
 \[
 T(x^i) K(\chi; x^i,p)=[K(\widetilde{\chi};x^i,p)+
  K_{\boldsymbol \mu,{\bold u}}(x^i) \Gamma^{-1} K^{{\bold x},
  \boldsymbol \lambda}(p)] T(p).
 \]
 We use this identity both in the form indicated and with $x^j$ in place 
 of $x^i$ to convert \eqref{Sside1} to
 \begin{align} \notag
 [S(p) \beta^{-1} u^{\times \prime}(p)]_i = &
 [K(\widetilde{\chi};x^i,p) +K_{\boldsymbol \mu,{\bold u}}(x^i) \Gamma^{-1} 
 K^{{\bold x},\boldsymbol \lambda}(p)] T(p) - \\
 \notag
  - K_{\boldsymbol \mu,{\bold u}}(x^i) \cdot & \underset{i^\prime}{\text{diag.}}
 \{(\xi_1 \lambda_1(p) + \xi_2 \lambda_2(p) - \xi_1 \mu^{i^\prime}_1 - 
 \xi_2 \mu^{i^\prime}_2)^{-1}I_{s_{i^\prime}}\} \cdot \\
 \label{Sside2}
  \cdot \Gamma^{-1} \cdot  \sum_{j=1}^m (\xi_1 c_{j1} + \xi_2 c_{j2}) &
 K^{{\bold x},\boldsymbol \lambda}(x^j) [K(\widetilde{\chi};x^j,p) 
 +K_{\boldsymbol \mu,{\bold u}}(x^j) \Gamma^{-1} 
 K^{{\bold x},\boldsymbol \lambda}(p)] T(p).
 \end{align}
 
 Next we use the general identity (see also \cite{AV})
 \begin{gather}
 \notag
 \sum_{j=1}^m (\xi_1 c_{j1} + \xi_2 c_{j2}) K(\widetilde{\chi};p,x^j) 
 K(\widetilde{\chi};x^j,q) = \\
 \label{collection}
 =(\xi_1 \lambda_1(q) + \xi_2 \lambda_2(q) - \xi_1 \lambda_1(p) - \xi_2 
 \lambda_2(p)) K(\widetilde{\chi}; p,q)
 \end{gather}
 which is valid for all distinct points $p,q$ in $X$ which are disjoint 
 from $x^1, \dots, x^m$.
 To prove this ``collection formula'' \eqref{collection}, consider each 
 side as a function of $p$ with $q$ fixed. Since $h^0(\widetilde{\chi} 
 \otimes \Delta) = 0$, it suffices to show that the local principal part
 in the Laurent series expansion at each pole of each side matches with 
 the local principal part of the other side.  One can check that the 
 only possible poles are all simple and occur at $x^1, \dots, x^m$ with 
 residue of each side at $x^i$ equal to the common value $-(\xi_1c_{i1} + 
 \xi_2 c_{i2}) dt^i(x^i) K(\widetilde{\chi}; x^i,q)$.
 Immediate consequences of the identity \eqref{collection} which are 
 important for our context here are:
 \begin{align} \notag
 \sum_{j=1}^m (\xi_1 c_{j1} + \xi_2 c_{j2}) & K^{{\bold x}, 
 \boldsymbol \lambda}(\lambda^i) K(\widetilde{\chi};x^j,p) = \\
 \label{consequence1}
 = & \underset{i^\prime}{\text{diag.}} \{ (\xi_1 \lambda_1(p) + \xi_2 
 \lambda_2(p) - \xi_1 \lambda^{i^\prime}_1 - \xi_2 
 \lambda^{i^\prime}_2)I_{t_{i^\prime}}\} \cdot 
 K^{{\bold x},\boldsymbol\lambda}(p),
 \end{align}
 and, if $(i^\prime, j^\prime)$ is a pair of indices for which 
 $\lambda^{i^\prime} \ne \mu^{j^\prime}$ then the $(i^\prime, 
 j^\prime)$-matrix entry of 
 $\sum_{j=1}^m (\xi_1 c_{j1} + \xi_2 c_{j2}) 
 K^{{\bold x}, \boldsymbol \lambda}(x^j) K_{\boldsymbol \mu,{\bold u}}(x^j)$
 is given by
 \begin{gather}  \notag
 [ \sum_{j=1}^m (\xi_1 c_{j1} + \xi_2 c_{j2}) 
 K^{{\bold x}, \boldsymbol\lambda}(x^j) 
 K_{\boldsymbol \mu,{\bold u}}(x^j)]_{i^\prime, j^\prime} 
 = \\ \notag
 = (\xi_1 \mu^{j^\prime}_1 + \xi_2 \mu^{j^\prime}_2 - \xi_1 
 \lambda^{i^\prime}_1 - \xi_2 \lambda^{i^\prime}_2) {\bold x}_{i^\prime} 
 K(\widetilde{\chi}; \lambda^{i^\prime}, \mu^{j^\prime}) {\bold u}_{j^\prime}
 \\ \notag
  = - (\xi_1 \mu^{j^\prime}_1 + \xi_2 \mu^{j^\prime}_2 - \xi_1 
 \lambda^{i^\prime}_1 - \xi_2 \lambda^{i^\prime}_2) \Gamma_{i^\prime, 
 j^\prime}.
 \end{gather}
 Hence in matrix form we have
 \begin{gather} \notag
   \sum_{j=1}^m (\xi_1 c_{j1} + \xi_2 c_{j2}) 
 K^{{\bold x}, \boldsymbol \lambda}(x^j) K_{\boldsymbol\mu,{\bold u}}(x^j) =
 \\ \label{consequence2}
 = \underset{i^\prime}{\text{diag.}} \{(\xi_1 \lambda^{i^\prime}_1
 + \xi_2 \lambda^{i^\prime}_2)I_{t_{i^\prime}} \} \cdot \Gamma - \Gamma 
 \cdot \underset{j^\prime}{\text{diag.}}\{(\xi_1 \mu^{j^\prime}_1 + \xi_2 
 \mu^{j^\prime}_2) I_{s_{j^\prime}} \}.
 \end{gather}
 
 In the case where $p=q$, the collection 
 formula \eqref{collection} takes the limiting form
 \begin{gather} \notag
 \sum_{j=1}^m (\xi_1 c_{j1} + \xi_2 c_{j2}) K(\widetilde{\chi};p,x^j) 
 K(\widetilde{\chi};x^j,p) = \\
 \label{degencollection}
 -(\xi_1 \lambda^\prime_1(p) + \xi_2 \lambda^\prime_2(p)) \ dt(p)
 \end{gather}
 where ${}^\prime = \dfrac{d}{dt}$ where $t$ is a local coordinate 
 centered at $p$.  An application of this degenerate collection formula
 \eqref{degencollection} 
 gives, for $\lambda^{i^\prime} = 
 \mu^{j^\prime} =: \xi^{i^\prime j^\prime}$,
 \begin{gather} \notag   
  [ \sum_{j=1}^m (\xi_1 c_{j1} + \xi_2 c_{j2}) 
 K^{{\bold x}, \boldsymbol \lambda}(x^j) 
 K_{\boldsymbol \mu,{\bold u}}(x^j)]_{i^\prime, j^\prime} 
 = \\ \notag
 = -(\xi_1 \lambda^\prime_1(\xi^{i^\prime j^\prime}) + \xi_2 
 \lambda^\prime_2(\xi^{i^\prime j^\prime}) )\ dt(\xi^{i^\prime j^\prime})
 \ {\bold x}^T_{i^\prime} {\bold u}_{j^\prime} = 0
 \end{gather}
 where we used the compatibility condition \eqref{comp} for the last 
 step.  We conclude that \eqref{consequence2} continues to be valid even 
 in the case where $\lambda^{i^\prime} = \mu^{j^\prime}$.
 
 Making the substitutions \eqref{consequence1} and \eqref{consequence2} 
 in \eqref{Sside2}, we obtain
 \begin{gather}  \notag
 [S(p) \beta^{-1} u^{\times \prime}(p)]_i =  
 [K(\widetilde{\chi};x^i,p)+
  K_{\boldsymbol \mu,{\bold u}}(x^i) \Gamma^{-1} 
  K^{{\bold x},\boldsymbol \lambda}(p)] T(p) -
  \\ \notag
  - K_{\boldsymbol\mu, {\bold u}}(x^i) \cdot  
  \underset{i'}{\text{diag.}}
 \{( \xi_1 \lambda_1(p)+\xi_2 \lambda_2(p) -  \xi_1 \mu^{i^\prime}_1
  -\xi_2 \mu^{i^\prime}_2)^{-1} I_{s_{i^\prime}} \} \cdot  
  \\  \notag
   \cdot 
  \Gamma^{-1} \cdot  \left( \underset{i^\prime}{\text{diag.}}\{(\xi_1 
  \lambda_1(p) + \xi_2\lambda_2(p) - \xi_1 \lambda^{i^\prime}_1 -\xi_2 
  \lambda^{i^\prime}_2) I_{t_{i^\prime}}\} 
  K^{{\bold x}, \boldsymbol \lambda}(p) + 
  \right.
  \\ \notag
   + \left[ \underset{i^\prime}{\text{diag.}}\{ (\xi_1 
  \lambda^{i^\prime}_1 + \xi_2 \lambda^{i^\prime}_2) I_{t_{i^\prime}} \}
  \Gamma \right.  - \Gamma \underset{j^\prime}{\text{diag.}} \{ ( \xi_1 
  \mu^{j^\prime}_1 + \xi_2 \mu^{j^\prime}_2) I_{s_{j^\prime}}\}
  \left. \right] 
  \Gamma^{-1} K^{{\bold x}, \boldsymbol\lambda}(p) \left. \right) =
  \\ \notag
  =  K(\widetilde{\chi};x^i,p) T(p)  + K_{\boldsymbol\mu, {\bold u}}(x^i) 
  \Gamma^{-1} K^{{\bold x}, \boldsymbol\lambda}(p) T(p) -
  \\ \notag
   - K_{\boldsymbol \mu, {\bold u}}(x^i) \cdot  
   \underset{i^\prime}{\text{diag.}}
  \{ (\xi_1 \lambda_1(p) + \xi_2 \lambda_2(p) - \xi_1 \mu^{i^\prime}_1 - 
  \xi_2 \mu^{i^\prime}_2)^{-1} I_{s_{i^\prime}} \} \cdot
  \\ \notag
   \cdot \left( (\xi_1 \lambda_1(p) + \xi_2 \lambda_2(p)) \Gamma^{-1} 
  K^{{\bold x}, \boldsymbol \lambda}(p)  \right. 
   - \Gamma^{-1} \underset{i^\prime}{\text{diag.}}
  \{(  \xi_1 \lambda^{i^\prime}_1 + \xi_2 \lambda^{i^\prime}_2) 
  I_{t_{i^\prime}} \} K^{{\bold x}, \boldsymbol \lambda}(p) + 
  \\ \notag
   + \Gamma^{-1} \underset{i^\prime}{\text{diag.}} 
   \{ (\xi_1 \lambda^{i^\prime}_1 + \xi_2 \lambda^{i^\prime}_2) 
   I_{t_{i^\prime}} \} 
  K^{{\bold x}, \boldsymbol \lambda}(p)  
  \left.
  - \underset{j^\prime}{\text{diag.}}
  \{(\xi_1 \mu^{j^\prime}_1 + \xi_2 \mu^{j^\prime}_2) I_{s_{j^\prime}} \} 
  \Gamma^{-1} K^{{\bold x}, \boldsymbol\lambda}(p) 
  \right) 
  \\ \notag
  =  K(\widetilde{\chi}; x^i,p) T(p)
  \end{gather}
  and the intertwining relation \eqref{intertwining} follows.  
\end{pf}

{\bf REMARK:}  Note that the results of this section in principle give a 
means of computing explicitly the unknown input bundle $\chi$ appearing 
in Theorem \ref{T:bvabsint} in Section \ref{S:absint}.  Indeed, given a 
data set $\omega$ for an (ABSINT) problem, we can convert it to a data 
set $\omega_0$ for a (CONINT) problem, as explained in the discussion 
preceding the statement of Theorem \ref{T:comparison}.  In the context of 
the (CONINT) problem, Theorem \ref{T:conint} solves the problem of 
identifying the input bundle.  Namely, the input bundle $E$ is defined to 
be $E=\ker (z_1\sigma_2 - z_2 \sigma_1 + \gamma)$ where $\sigma_1, 
\sigma_2 $ and $\widetilde{\gamma}$ are given by \eqref{pencilcoef} and 
$\gamma$ is given by \eqref{gamma}.  The vector bundle $\chi$ is then 
determined up to biholomorphic equivalence by the condition
\[
\chi \cong E \otimes {\cal O}(1) \otimes \Delta.
\]
Moreover, as explained in \cite{hip}, it is possible to construct a 
matrix of normalized sections $u^\times$ for $E$ from working with the 
minors of the matrix pencil $z_1 \sigma_2 - z_2 \sigma_1 + \gamma$.  Such 
a matrix of normalized sections $u^\times$ in turn implements concretely 
the biholomorphic equivalence between $\chi$ and $E \otimes {\cal O}(1) 
\otimes \Delta$.

\end{document}